\documentclass[a4paper]{article}

\usepackage{amssymb}
\usepackage{amsmath}
\usepackage{amsthm}
\usepackage{bm}
\usepackage{color}
\usepackage{dsfont}
\usepackage{epic}
\usepackage{eepic}
\usepackage{enumerate}
\usepackage[OT2,OT1]{fontenc}
\usepackage{graphicx}
\usepackage{ifthen}
\usepackage{pifont}
\usepackage{textcomp}
\usepackage{wasysym}
\usepackage{todonotes}

\newcommand{\mc}[1]{{\mathcal{#1}}}

\newcommand{\bb}[1]{{\mathbb{#1}}}

\makeatletter
\newcommand\appendixsection{\@startsection {section}{1}{\z@}
	{-3.5ex \@plus -1ex \@minus -.2ex}{2.3ex \@plus.2ex}
	{\normalfont\Large\bfseries\hspace*{-15pt}Appendix\ }}
\makeatother

\numberwithin{equation}{section}
\swapnumbers
\theoremstyle{plain}
\newtheorem{lemma}{Lemma}[section]
\newtheorem{proposition}[lemma]{Proposition}
\newtheorem{theorem}[lemma]{Theorem}
\newtheorem{corollary}[lemma]{Corollary}
\theoremstyle{definition}
\newtheorem{definition}[lemma]{Definition}

\newtheorem{remark}[lemma]{Remark}

\newcommand\cyr{%
	\renewcommand\rmdefault{wncyr}%
	\renewcommand\sfdefault{wncyss}%
	\renewcommand\encodingdefault{OT2}%
	\normalfont
	\selectfont}
\DeclareTextFontCommand{\textcyr}{\cyr}

\DeclareMathOperator{\diag}{diag}
\DeclareMathOperator{\dom}{dom}
\DeclareMathOperator{\loc}{loc}

\DeclareMathOperator{\mul}{mul}
\DeclareMathOperator{\ran}{ran}
\DeclareMathOperator{\mult}{mult}
\DeclareMathOperator{\rank}{rank}

\DeclareMathOperator{\spn}{span}

\DeclareMathOperator{\tr}{tr}
\DeclareMathOperator{\Cls}{cls}

\DeclareMathOperator{\IM}{Im}

\renewcommand{\Im}{\IM}

\newlength{\maxlabwidth}

\begin{document}
	{\Large\bf
		\begin{flushleft}
			Local spectral multiplicity of selfadjoint couplings with 
			general interface conditions
		\end{flushleft}
	}
	\vspace*{3mm}
	\begin{center}
		{\sc Sergey Simonov, Harald Woracek}
	\end{center}
	
	\begin{abstract}
		\noindent
		We consider selfadjoint operators obtained by pasting a finite number of boundary relations with one-dimensional
		boundary space. A typical example of such an operator is the Schr\"odinger operator on a star-graph with a finite
		number of finite or infinite edges and an interface condition at the common vertex. A wide class of
		``selfadjoint'' interface conditions, subject to  an assumption which is generically satisfied, is considered. 
		We determine the spectral multiplicity function on the singular spectrum (continuous as well as point) 
		in terms of the spectral data of decoupled operators.
		
	\end{abstract}
	\begin{flushleft}
		{\small
			{\bf AMS MSC 2010:} 34\,B\,45, 34\,B\,20, 34\,L\,40, 47\,E\,05, 47\,J\,10 
		}
	\end{flushleft}

%%%%%%%%%%%%%%%%%%%%%%%%%%%%%%%%%%%%%%%%%%%%%%%%%%%%%%%%%%%%%%%%%%%%%%%%%%%%%%%%%%%%%%%%%%%%%%%%%%%

\section{Introduction}\label{section introduction}

In this paper we analyze the singular spectrum of a selfadjoint operator built by gluing together a finite number of
selfadjoint operators with simple spectrum by means of interface conditions. We realize the operators with help of
boundary triplets, and understand interface conditions as linear dependencies among boundary values. An archetypical
example for such a glued together operator is a Schr\"odinger operator on a star-graph with an interface condition at the inner
vertex. 

The reader who is not familiar with the language of boundary triplets and couplings of such (e.g.\ 
\cite{Behrndt-Hassi-deSnoo-2020,Derkach-Hassi-Malamud-deSnoo-2009}), is advised to pause for a moment, and 
before proceeding here read through Section~1.1 below in order to get an intuition; there we elaborate in detail the above 
mentioned example. We also point out that all necessary notions and results from the literature are recalled in Section~2.

Assume we have selfadjoint operators $L_l$, $l=1,\ldots,n$, in Hilbert spaces $H_l$, $l=1,\ldots,n$, that emerge from boundary
relations with one-dimensional boundary value space (cf.\ Section~2.3). Denote by $L_0$ the orthogonal coupling, i.e., the 
diagonal operator, $L_0:=\prod_{l=1}^n L_l$ acting in $H:=\prod_{l=1}^n H_l$. The spectrum of $L_0$ and its multiplicity 
is easily understood: letting $N_l$ and $N_0$ be the respective spectral multiplicity functions  of
$L_l$ and $L_0$, it holds that $N_0=\sum_{l=1}^n N_l$. 
Here -- and always -- we tacitly understand that any relation between spectral multiplicity functions should be valid only after
making an appropriate choice of representants for each of them. 

The orthogonal coupling can be seen as gluing together the single operators without allowing any interaction between them.
Much more interesting is what happens when the single operators do
influence each other after gluing. Assume we have an interface condition, formulated as a linear dependence among boundary values 
and described by matrices $A,B$ (cf.\ Section~2.4), such that the corresponding operator $L_{A,B}$ is selfadjoint, 
and let $N_{A,B}$ be the spectral multiplicity function of $L_{A,B}$. The basic question is: 
\begin{center}
	{\it
	How to compute $N_{A,B}$ from $N_l$, $l=1,\ldots,n$ ?
	}
\end{center}
By a simple dimension argument, it always holds that $N_{A,B}(x)\leq n$. 
Further, the Kato-Rosenblum theorem fully settles the question on the absolutely continuous part of the spectrum: 
there we have $N_{A,B}(x)=N_0(x)$ since $L_{A,B}$ and $L_0$ are finite rank perturbations of each other in the resolvent sense. 
Contrasting this, the singular spectrum (eigenvalues as well as singular continuous part) may behave wildly, and much less is 
known. 

A classical result for the case that $n=2$ is Kac's Theorem \cite{Kac-1962,Kac-1963}, which says in the present language 
that for a particular interface condition, namely the standard condition, the multiplicity of the singular spectrum does not
exceed $1$. His proof proceeds via an analysis of the Cauchy transforms of the involved spectral measures and the
Titchmarsh-Kodaira formula. Different proofs are given in \cite{Gilbert-1998} (by using subordinacy theory) and in 
\cite{Simon-2005} (by reducing to the Aronszajn-Donoghue Theorem). Generalizations of the theorem of Kac are given in 
\cite{Simonov-Woracek-2014} and in \cite{Malamud-2019}. In the first reference, we allow arbitrary $n$ but still prescribe the
standard interface condition. The second reference goes into another direction. There still $n=2$ but the boundary value
spaces are allowed to have higher dimensions and a certain class of interface conditions is permitted. 

In the present paper we allow arbitrary $n$ and consider a fairly rich class of interface conditions defined by an algebraic 
property (cf.\ Section~3). 
This property expresses, at least in some sense, that all single operators influence each other and no splitting into
independent blocks happens, though one has to be careful with this intuition, it is only very rough. The previously considered 
standard condition belongs to that class. A striking difference is that interface conditions of the presently considered class 
can give rise to perturbations of any rank up to $n$, while the standard condition always yields a rank one perturbation. 
Our main result says that, letting $r$ be the rank of the perturbation, on the singular spectrum the relation
\begin{equation}
\label{aaa}
	N_{A,B}(x)\ 
	\begin{cases}
		\ =N_0(x)-r &\text{if } N_0(x)\geq r,
		\\[2mm]
		\ \leq r-N_0(x) &\text{if } N_0(x)<r,
	\end{cases}
\end{equation}
holds. 
A formulation in terms of spectral measures, and without reference to a particular choice of representants of multiplicity
functions, is given in Theorem~\ref{theorem main}. 
The proof of the theorem is carried out by an in depth analysis of Cauchy integrals and the matrix measure in the
Titchmarsh-Kodaira formula. We exploit algebraic properties of the considered class of interface conditions to obtain the
rank of the derivative of that matrix measure w.r.t.\ its trace measure, and this leads to \eqref{aaa}.

Other approaches to the above emphasised basic question might proceed via the already mentioned work of M.Malamud
\cite{Malamud-2019}, or via a generalization of Aronzsajn-Donoghue's theorem given by C.Liaw and S.Treil in 
\cite{Liaw-Treil-2018}. To the best of our knowledge, no  such results have been obtained so far using these approaches.

The present paper is organized as follows. 
In Section~\ref{section preliminaries} we introduce objects and tools needed to formulate and
prove our result. In particular, these include boundary relations, pasting of such with selfadjoint interface conditions,
matrix-valued Weyl functions and corresponding measures. In Section~\ref{section interface condition} we discuss in detail the
class of selfadjoint interface conditions that we consider, namely, their description in terms of matrices $A$, $B$ and the
additional assumption that we make about them. In Section \ref{section formulation} we give the statement of the main result,
Theorem~\ref{theorem main}. Before that we formulate the result separately for the case of point spectrum, 
Theorem~\ref{theorem point}, since this can be shown under slightly weaker assumptions. 
In Section \ref{section point} we prove Theorem~\ref{theorem point}. In Section~\ref{section many layers}
we prove the part of Theorem~\ref{theorem main} concerning the case where many layers of the spectrum ``overlap''. 
In Section~\ref{section few layers} we prove the remaining part of Theorem~\ref{theorem main}.

\subsection{The Schr\"odinger operator on a star-graph}

Let  us discuss, as an example, the Schr\"odinger operator on a metric star-graph. In fact, this example can 
serve as a model for the general case. 
We denote the edges of the graph by $E_1,\ldots,E_n$ and associate them with intervals $[0,e_l)$, where the
endpoint $0$ corresponds to the inner vertex and $e_l$ can be finite or infinite. Assume we are given data:
\begin{enumerate}[(1)]
\item Real-valued potentials $q_l\in L_{1,\loc}([0,e_l))$ for $l=1,\ldots,n$.
\item Boundary conditions at $e_l$ for those $l=1,\ldots,n$, for which $q_l$ is regular or is in the limit point case at $e_l$.
\end{enumerate}
For $l=1,\ldots,n$ let $H_l$ be the Hilbert space $L_2(0,e_l)$, and let $L_l$ be the selfadjoint Schr\"odinger operator with 
Dirichlet boundary conditions: 
$$
	L_lu:=-u''+q_lu,
$$
$$
	\begin{aligned}
		\dom L_l := \Bigg\{
		&u\in L_2(0,e_l):
		u,u'\text{ are absolutely continuous},\\
		& -u''+q_lu\in L_2(0,e_l),\quad u(0)=0,\\
		& u\text{ satisfies the boundary condition at $e_l$ (if present)}
		\Bigg\}\ .
	\end{aligned}
$$
Now assume we have an interface condition at the inner vertex written in the form 
\begin{equation}\label{i1}
	A
	\begin{pmatrix}
		u_1(0)\\
		\vdots\\
		u_n(0)
	\end{pmatrix}
	+B
	\begin{pmatrix}
		u_1'(0)\\
		\vdots\\
		u_n'(0)
	\end{pmatrix}=0,
\end{equation}
where $A$ and $B$ are $n\times n$ matrices such that
\begin{equation}\label{i2}
	AB^*=BA^*,\quad{\rm rank\,}(A,B)=n.
\end{equation}
Here $(A,B)$ denotes the $n\times2n$ matrix which has $A$ as its first $n$ columns and $B$ as its last $n$ columns. 
The operator $L_{A,B}$ is defined in the Hilbert space $H:=\prod_{l=1}^n L_2(0,e_l)$ and acts by the rule
\begin{equation}\label{i3}
	L_{A,B}
	\begin{pmatrix}
		u_1\\
		\vdots\\
		u_n
	\end{pmatrix}
	:=
	\begin{pmatrix}
		-u_1''\\
		\vdots\\
		-u_n''
	\end{pmatrix}
	+
	\begin{pmatrix}
		q_1u_1\\
		\vdots\\
		q_nu_n
	\end{pmatrix}
\end{equation}
on the domain
\begin{equation}\label{i4}
	\begin{aligned}
		\dom L_{A,B} := \Bigg\{
		&(u_1,\dots,u_n) \in\prod_{l=1}^n L_2(0,e_l):\quad \forall l\in\{1,\ldots,n\}\\
		& u_l,u_l'\text{ are absolutely continuous}, -u_l''+q_lu_l\in L_2(0,e_l),\\[2mm]
		& u_l\text{ satisfies the boundary condition at $e_l$ (if present),}\\[-1mm]
		& u_1,\dots,u_n\text{ satisfy the interface condition \eqref{i1}}
		\Bigg\}\ .
	\end{aligned}
\end{equation}
Since the matrices $A$ and $B$ satisfy \eqref{i2}, the operator $L_{A,B}$ is selfadjoint \cite{Kostrykin-Schrader-1999}. 
Obviously, this correspondence between matrices and operators is not one-to-one: one can multiply $A$ and $B$ simultaneously
from the left by any invertible matrix, and this defines the same interface condition and the same operator.

The orthogonal coupling $L_0:=\prod_{l=1}^nL_l$ corresponds to the matrices $A_0=I$, $B_0=0$: in the above notation 
$L_0=L_{I,0}$. The standard interface condition corresponds to the matrices
\begin{equation}\label{standard ic}
	A_{st}=
	\begin{pmatrix}
		1&\hdots&0&-1\\
		\vdots&\ddots&\vdots&\vdots\\
		0&\hdots&1&-1\\
		0&\hdots&0&0
	\end{pmatrix},
	\quad
	B_{st}=
	\begin{pmatrix}
		0&\hdots&0&0\\
		\vdots&\ddots&\vdots&\vdots\\
		0&\hdots&0&0\\
		1&\hdots&1&1
	\end{pmatrix}
\end{equation}
where the first $n-1$ lines express continuity of the solution at the vertex, and the last line corresponds to the Kirchhoff
condition that the sum of derivatives vanishes.
The class of interface conditions we consider in the present paper is given by those matrices $A,B$ subject to \eqref{i2} 
which satisfy in addition the following assumption: each set of $\rank B$ many different columns of $B$ is lineary independent.
Under this assumption the rank of the difference of resolvents of $L_{A,B}$ and $L_0$ equals ${\rm rank\,}B$, and 
\eqref{aaa} holds.

%%%%%%%%%%%%%%%%%%%%%%%%%%%%%%%%%%%%%%%%%%%%%%%%%%%%%%%%%%%%%%%%%%%%%%%%%%%%%%%%%%%%%%%%%%%%%%%%%%%

\subsection*{Acknowledgements}
The first author was supported by the RScF 20-11-20032 grant (Sections \ref{section introduction}--\ref{section formulation}), 
by the RFBR 19-01-00565A and by the RFBR 19-01-00657A grants (Sections \ref{section point}--\ref{section few layers}). 
He appreciates the hospitality of TU Wien where the most part of this work was done.

The second author was supported by the joint project I~4600 of 
the Austrian Science Fund (FWF) and the Russian foundation of basic research (RFBR).

%%%%%%%%%%%%%%%%%%%%%%%%%%%%%%%%%%%%%%%%%%%%%%%%%%%%%%%%%%%%%%%%%%%%%%%%%%%%%%%%%%%%%%%%%%%%%%%%%%%

\section{Preliminaries}\label{section preliminaries}

\subsection{Boundary behavior of Herglotz functions}

Recall the notion of matrix valued Herglotz functions (often also called Nevanlinna functions).

\begin{definition}
	An analytic function $M:\bb C\setminus\bb R\to\bb C^{n\times n}$ is called a (\emph{$n\times n$-matrix valued}) 
	\emph{Herglotz function}, if
	\begin{enumerate}[(H1)]
	\item $M(\overline z)=M(z)^*$, $z\in\bb C\setminus\bb R$.
	\item For each $z\in\bb C^+$, the matrix $\Im M(z):=\frac 1{2i}(M(z)-M(z)^*)$ is positive semidefinite.
	\end{enumerate}
\end{definition}

\noindent
Any Herglotz function $M$ admits an integral representation. 
Namely, there exists a finite positive $n\!\times\!n$-matrix valued Borel measure $\Omega$ (which means that $\Omega(\Delta)$ is a positive semidefinite matrix for every Borel set $\Delta$), a selfadjoint matrix $a$, and a positive 
semidefinite matrix $b$, such that
\begin{equation}
\label{equ 20}
	M(z)=a+bz+\int_{\bb R}\frac{1+xz}{x-z}\,d\Omega(x),\quad z\in\bb C\setminus\bb R
	.
\end{equation}
For the scalar case, this goes back to \cite{Herglotz-1911}, for the matrix valued case see, e.g., \cite[Theorem 5.4]{Gesztesy-Tsekanovskii-2000}.

We use several known facts about the boundary behavior of Herglotz functions which relate normal or nontangential boundary limits to the measure $\Omega$ in \eqref{equ 20}. The key notion in this context is the symmetric derivative of one measure
relative to another. If $\sigma$ is a positive Borel measure and $\nu$ is a positive Borel measure or a complex Borel measure 
absolutely continuous w.r.t. $\sigma$, then we define the \emph{symmetric derivative} $\frac{d\nu}{d\sigma}(x)$ at a point $x\in\bb R$ as the limit
\[
	\frac{d\nu}{d\sigma}(x):=\lim\limits_{\varepsilon\to0+}
	\frac{\nu([x-\varepsilon,x+\varepsilon])}{\sigma([x-\varepsilon,x+\varepsilon])}
	,
\]
whenever it exists in $[0,\infty]$, or in $\bb C$, respectively. 

\begin{proposition}[\cite{diBenedetto-2002}]\label{prop DiB02}
	There exists a Borel set $X\subseteq\bb R$ with $\sigma(\bb R\setminus X)=0$, such that the symmetric derivative exists for all $x\in X$ and the function $\frac{d\nu}{d\sigma}$ is measurable on $X$. 
\end{proposition}

By the de la Valle\'e--Poussin theorem \cite{saks-1964,diBenedetto-2002} the function $\frac{d\nu}{d\sigma}$ is a 
Radon--Nikodym derivative of $\nu$
with respect to $\sigma$. We formulate two corollaries of the de la Valle\'e--Poussin theorem which will be of particular
convenience to us in what follows. An explicit proof can be found in \cite{Simonov-Woracek-2014}.

The first corollary concerns properties of sets. A set $X\subseteq\bb R$ is called $\nu$-{\it zero}, if
there exists a Borel set $X'\supseteq X$ such that $\nu(X')=0$; a set is called $\nu$-{\it full}, if its complement is
$\nu$-zero. For a Borel set $X$ the measure $\mathds{1}_X\cdot \nu$ is defined as
$(\mathds{1}_X\cdot\nu)(\Delta)=\nu(X\cap\Delta)$ on Borel sets $\Delta$.
\begin{corollary}\label{prop Vallee-Poussin sets}
	Let $\nu$ and $\sigma$ be positive Borel measures, and let $X\subseteq\bb R$.
	\begin{enumerate}[(1)]
		\item If $\frac{d\nu}{d\sigma}(x)=0$ for all $x\in X$, then $X$ is $\nu$-zero.
		\item If the set $X$ is $\nu$-zero, then $\frac{d\nu}{d\sigma}(x)=0$ for $\sigma$-a.a. $x\in X$.
		\item If $X$ is a Borel set and $\frac{d\nu}{d\sigma}(x)\in[0,\infty)$ for all $x\in X$, then $\mathds{1}_X\cdot\nu\ll\sigma$.
		\item If $X$ is a Borel set and $\frac{d\nu}{d\sigma}(x)\in(0,\infty)$ for all $x\in X$, then $\mathds{1}_X\cdot\nu\sim\mathds{1}_X\cdot\sigma$.
	\end{enumerate}
\end{corollary}
The second corollary concerns properties of the symmetric derivative.
\begin{corollary}\label{prop Vallee-Poussin measures}
	Let $\nu$ and $\sigma$ be positive Borel measures on $\bb R$.	Let $\nu=\nu_{ac}+\nu_s$ and $\sigma=\sigma_{ac}+\sigma_s$
	be the Lebesgue decompositions of $\nu$ with respect to $\sigma$ and of $\sigma$ with respect to $\nu$, respectively.  Then
	\begin{enumerate}[(1)]
		\item $\frac{d\nu}{d\sigma}(x)\in[0,\infty)$, $\sigma$-a.e.
		\item $\frac{d\nu}{d\sigma}(x)\in(0,\infty]$, $\nu$-a.e.
		\item $\frac{d\nu}{d\sigma}(x)\in(0,\infty)$, $\nu_{ac}$-a.e. and $\sigma_{ac}$-a.e.
		\item $\frac{d\nu}{d\sigma}(x)=\infty$, $\nu_s$-a.e.
		\item $\frac{d\nu}{d\sigma}(x)=0$, $\sigma_s$-a.e.
	\end{enumerate}
\end{corollary}
The following statements concern the relationship between boundary behavior of Herglotz functions and symmetric derivatives of the measures associated with them.

\begin{proposition}[\cite{Pearson-1988}, \cite{Kac-1962}, \cite{Poltoratski-1994}, \cite{Poltoratski-2009}]
\label{proposition 20}
	\phantom{}
	\begin{enumerate}[(1)]
	\item Let $\nu$ and $\sigma$ be finite positive Borel measures and $x\in\bb R$. Assume that $\frac{d\nu}{d\sigma}(x)$ exists in $[0,\infty)$ 
		and $\frac{d\sigma}{d\lambda}(x)$ exists in $(0,\infty]$. Let $m_\nu$ and $m_\sigma$ be two Herglotz functions 
		with the measures $\nu$ and $\sigma$, respectively, in their integral representations \eqref{equ 20}. Then 
		\[
			\lim_{z\downarrow x}\frac{\Im m_\nu(z)}{\Im m_\sigma(z)}=\frac{d\nu}{d\sigma}(x)
			.
		\]
	\item Let $\nu$ be a finite positive Borel measure, let $m_\nu$ be a Herglotz function with 
		the measure $\nu$ in its integral representation, and let $x\in\bb R$. If $\frac{d\nu}{d\lambda}(x)=\infty$, then 
		\[
			\lim_{z\downarrow x}\Im m_\nu(z)=\infty
			.
		\]
	\item Let $\nu$ and $\sigma$ be finite positive Borel measures with $\nu\ll\sigma$, and let $\sigma_s$ be the singular part of
		$\sigma$ w.r.t.\ $\lambda$. Again let $m_\nu$ and $m_\sigma$ be two Herglotz functions with $\nu$ or $\sigma$,
		respectively, in their integral representations. Then 
		\[
			\lim_{z\downarrow x}\frac{m_\nu(z)}{m_\sigma(z)}=\frac{d\nu}{d\sigma}(x)
			\text{ for $\sigma_s$-a.a. }x\in\bb R
			.
		\]
	\end{enumerate}
\end{proposition}

\noindent
Note that in (1) we can use $d\sigma:=\frac{d\lambda}{1+x^2}$, which implies
\[
	\lim_{z\downarrow x}\Im m_\nu(z)=\pi(1+x^2)\frac{d\nu}{d\lambda}(x)
\]
whenever $\frac{d\nu}{d\lambda}(x)$ exists and is finite. 

For a more detailed compilation and references about symmetric derivatives and boundary values we refer the reader to 
\cite[\S2.3 and \S2.4]{Simonov-Woracek-2014}.

Item (1) of the above theorem has an obvious extension to matrix valued functions and measures. 

\begin{lemma}\label{lemma matrix Herglotz}
\label{lemma 20}
	Let $M$ be a $n\times n$-matrix valued Herglotz function and let $\Omega$ be the measure in its integral representation 
	\eqref{equ 20}. Denote by $\rho$ the trace measure of $\Omega$, i.e., $\rho(\Delta):=\tr\Omega(\Delta)$ for every Borel
	set $\Delta$. Then $\Omega\ll\rho$, and for $\rho$-a.a.\ $x\in\bb R$ the symmetric derivative $\frac{d\Omega}{d\rho}(x)$
	exists and is related to $M$ by 
	\begin{equation}\label{p17}
		\lim_{z\downarrow x}\frac{\Im M(z)}{\Im\tr M(z)}=\frac{d\Omega}{d\rho}(x)
		.
	\end{equation}
\end{lemma}
\begin{proof}
	Let $x\in\bb R$, and assume that all symmetric derivatives at $x$ of positive and negative parts of the real and the imaginary 
	parts of entries of $\Omega$ w.r.t.\ $\rho$ exist and are finite. 
	This is fulfilled for $\rho$-a.a.\ $x\in\bb R$. 

	We have $M(\bar z)=M(z)^*$ and hence
	\begin{align*}
		\Im M(z)= &\, \frac 1{2i}\big(M(z)-M(z)^*\big)=\frac 1{2i}\big(M(z)-M(\bar z)\big)
		\\
		= &\, b\Im z+\int_{\bb R}\Im\Big(\frac{1+xz}{x-z}\Big)\,d\Omega(x)
		.
	\end{align*}
	Hence, using linearity of the integral in the measure, Proposition~\ref{proposition 20}, (1), and 
	Corollary \ref{prop Vallee-Poussin measures}, (1),(2), we obtain 
	\[
		\lim_{z\downarrow x}\frac{\Im M_{kl}(z)}{\Im\tr M(z)}=\frac{d\Omega_{kl}}{d\rho}(x)\text{ for $\rho$-a.a. }x\in\mathbb R
		,\quad k,l=1,\ldots n
		.
	\]
\end{proof}

\subsection{Spectral multiplicity}\label{section multiplicity}

Let $L$ be a selfadjoint operator in a Hilbert space $H$ and $E_L$ be its projection valued spectral measure. 
A subspace $G$ of $H$ is called generating, if 
\[
	\Cls\left(\bigcup\left\{ E_L(\Delta)G:\Delta\text{ Borel set}\right\}\right)=H
	.
\]
The minimum of dimensions of generating subspaces is called the (global) {\it spectral multiplicity} of operator $L$ and is denoted by $\mult L$. 

The operators that we deal with always have finite spectral multiplicity, hence we shall assume from now on that 
$\mult L<\infty$. There exist elements $g_1,\ldots,g_{\mult L}\in H$, such that their linear span is a generating subspace and the positive measures
\begin{equation}\label{scalar measure}
\nu_l(\Delta):=(E_L(\Delta)g_l,g_l),\quad l=1,\ldots,\mult L,
\end{equation}
are ordered in the sense of absolute continuity as 
\[
	\nu_{\mult L}\ll\cdots\ll\nu_1\sim E_L
	.
\]
An element $g_1$ occuring in a collection of elements with these properties is called an
{\it element of maximal type}. 

Fix a choice of such elements $g_1,\ldots,g_{\mult L}$. Then the operator $L$ is unitarily equivalent to the operator of multiplication by the independent variable in the space $\prod_{l=1}^{\mult L}L_2(\mathbb R,\nu_l)$. We choose Borel sets $Y_1,\ldots,Y_{\mult L}$ such that 
\[
	Y_{\mult L}\subseteq\ldots\subseteq Y_1,\quad \nu_l(\bb R\setminus Y_l)=0,\quad \mathds{1}_{Y_l}\cdot\nu_l\sim\mathds{1}_{Y_l}\cdot\nu_1
	,
\]
and define the (local) spectral multiplicity function of $L$ as the equivalence class (i.e., functions which are $E_L$-a.e. equal) of the function 
\begin{equation}\label{p2}
	N_L(x):=\#\big\{l:x\in Y_l\big\}
	.
\end{equation}
The need of considering an equivalence class of functions arises since the sets $Y_l$ are defined non-uniquely up to $\nu_1$-zero sets and thus the function $\#\{l:x\in Y_l\}$ can be changed on $\nu_1$-zero sets.

Intuitively, the sets $Y_l$ correspond to the ``layers'' of the spectrum, and hence indeed $N_L$ expresses spectral
multiplicity in a natural way. If $x$ is an eigenvalue of $L$, then $N_L(x)={\rm dim\,ker\,}(L-xI)$ is the usual multiplicity of an eigenvalue.

The spectral multiplicity function is a unitary invariant of the operator and does not depend on a choice of generating basis.

\begin{remark}
	We will also speak of the spectral multiplicity function of a selfadjoint linear relation $L$ in a Hilbert space $H$. 
	This notion is defined by simply ignoring the multivalued part (and doing so is natural, since one can think of the 
	multivalued part as an eigenspace for the eigenvalue $\infty$). To be precise, let $\mul L:=\{g\in H:(0;g)\in L\}$. 
	Then 
	\[
		L_{op}:= L\cap\big((\mul L)^\perp\times(\mul L)^\perp\big)
	\]
	is a selfadjoint operator (recall that we identify operators with their graphs) in the Hilbert space $(\mul L)^\perp$, and 
	$L=L_{op}\oplus(\{0\}\times\mul L)$. Now we define $N_L:=N_{L_{op}}$. 
\end{remark}

The following classical fact will be used below. 

\begin{lemma}\label{proposition orthogonal sum}
	Let $L_1$ and $L_2$ be selfadjoint relations in Hilbert spaces $H_1$ and $H_2$, respectively, with finite multiplicities. Set 
	$H:=H_1\oplus H_2$ and $L:=L_1\oplus L_2$. Let $\mu=(E_{L_1}f,f)\sim E_{L_1}$, $\nu=(E_{L_2}g,g)\sim E_{L_2}$ be scalar measures defined 
	by elements of maximal type $f$ for $L_1$ and $g$ for $L_2$ via \eqref{scalar measure}. Let
	$$
	\mu=\mu_{ac}+\mu_s,\quad\nu=\nu_{ac}+\nu_s
	$$
	be the Lebesgue decompositions of the measures $\mu$, $\nu$ with respect to each other: 
	$\mu_{ac}\ll\nu$, $\mu_s\perp\nu$; $\nu_{ac}\ll\mu$, $\nu_s\perp\mu$, $\mu_{ac}\sim\nu_{ac}$. 
	Then 
	$$
	\begin{array}{ll}
		N_L=N_{L_1}, &\quad \mu_s\text{-a.e.},
		\\
		N_L=N_{L_1}+N_{L_2}, &\quad \mu_{ac}\text{-a.e. and }\nu_{ac}\text{-a.e.},
		\\
		N_L=N_{L_2}, &\quad\nu_s\text{-a.e.}
	\end{array}
	$$	
\end{lemma}

\noindent
Consider the measure $\mu+\nu$ and the sets
$$
X_1:=\left\{ x\in\bb R:\frac{d\mu}{d(\mu+\nu)}(x)>0\right\},\quad X_2:=\left\{ x\in\bb R:\frac{d\nu}{d(\mu+\nu)}(x)>0\right\},
$$
 Then the sets $X_1\setminus X_2$, $X_1\cap X_2$ and $X_2\setminus X_1$ carry the measures $\mu_s$, $\mu_{ac}\sim\nu_{ac}$ and $\nu_s$, respectively. We can see that caution is needed when dealing with local spectral
multiplicities: the values of function $N_{L_2}$  have no meaning on the set $X_1\setminus X_2$ and can be changed arbitrarily
there, the same for $N_{L_1}$ on $X_2\setminus X_1$. So the statement $N_L=N_{L_1}+N_{L_2}$, looking natural, is in fact
wrong: the functions $N_{L_1}$ and $N_{L_2}$ are defined non-uniquely, each in its own sense.

The next lemma is a general and folklore-type result for the behavior of local spectral multiplicity of an operator under 
finite-dimensional perturbation generalized to the linear relations case. We do not know an explicit reference for it, and for the  reader's convenience we provide its proof.

\begin{lemma}\label{lemma unknown}
Let $L_1$ and $L_2$ be selfadjoint relations in the Hilbert space $H$ such that for ( some $\Leftrightarrow$ every) $\lambda\in\bb C\setminus\bb R$
$$
	\rank\big[(L_1-\lambda I)^{-1}-(L_2-\lambda I)^{-1}\big]=k.
$$
Let $E_1, E_2$ be projection valued spectral measures of their operator parts, $\mu=(E_{L_1}f,f)\sim E_{L_1}$, $\nu=(E_{L_2}g,g)\sim E_{L_2}$ be scalar measures defined 
by elements of maximal type $f$ for $L_1$ and $g$ for $L_2$ via \eqref{scalar measure} and $N_1,N_2$ be their local spectral multiplicity functions.
Let
$$
\mu=\mu_{ac}+\mu_s,\quad\nu=\nu_{ac}+\nu_s
$$
be Lebesgue decompositions of the measures $\mu$, $\nu$ with respect to each other: 
$\mu_{ac}\ll\nu$, $\mu_s\perp\nu$; $\nu_{ac}\ll\mu$, $\nu_s\perp\mu$, $\mu_{ac}\sim \nu_{ac}$. Then
\begin{enumerate}[\rm(1)]
	\item $N_1\leqslant k$, $\mu_s$-a.e.,
	\item $|N_1-N_2|\leqslant k$, $\mu_{ac}$-a.e. and $\nu_{ac}$-a.e.,
	\item $N_2\leqslant k$, $\nu_s$-a.e.
\end{enumerate}
\end{lemma}

\begin{proof}
Consider the symmetric linear relation $S:=L_1\cap L_2$ (recall once more that we identify operators with their graphs). It has
an orthogonal decomposition into a selfadjoint and a simple (i.e., completely nonselfadjoint) symmetric part. 
The first, a selfadjoint linear relation $L$, acts in the subspace 
$$
	H_L:=\bigcap_{\lambda\in\bb C\setminus\bb R}\ran(S-\lambda I)
	.
$$
The second, a simple symmetric operator $\widetilde S$, acts in the subspace
$$
	\widetilde H_S:=\Cls\Big(\bigcup_{\lambda\in\bb C\setminus\bb R}\ker(S^*-\lambda I)\Big)
	.
$$
The subspaces $H_L$ and $\widetilde H_S$ reduce the relations $S$, $L_1$, $L_2$, and $S^*$, see 
\cite[Lemma 3.4.2]{Behrndt-Hassi-deSnoo-2020}, and $S=L\oplus \widetilde S$. Thus $L_1$ and $L_2$ have orthogonal 
decompositions
$$
	L_1=L\oplus \widetilde L_1,\quad L_2=L\oplus \widetilde L_2,
$$
where the linear relations $\widetilde L_1$ and $\widetilde L_2$ are selfadjoint extensions of $\widetilde S$. 

Consider $\lambda\in\bb C\setminus\bb R$ and the subspace
$$
	H_{\lambda}:=\ker\big[(L_1-\lambda I)^{-1}-(L_2-\lambda I)^{-1}\big].
$$
We have: for every $v\in H_{\lambda}$
$$
	u:=(L_1-\lambda I)^{-1}v=(L_2-\lambda I)^{-1}v,
$$
hence $(u;v)\in(L_1-\lambda I)\cap(L_2-\lambda I)=S-\lambda I$ and $v\in\ran(S-\lambda I)$. The converse is also true, and therefore
$$
	H_{\lambda}=\ran(S-\lambda I).
$$
By assumption the subspace $H_{\lambda}$ has codimension $k$, hence the deficiency index of 
$S$ is $(k,k)$. The deficiency index of $\widetilde S$ coincides with that of $S$ and hence is $(k,k)$. Then spectral multiplicities of both relations $\widetilde L_1$ and $\widetilde L_2$ do not exceed $k$, because defect subspaces of simple symmetric operators are generating subspaces for the operator parts of their selfadjoint extensions. 

The rest follows from Lemma \ref{proposition orthogonal sum}. One should write out the ``triple'' Lebesgue decomposition for
scalar spectral measures of maximal type of $L$, $\widetilde L_1$ and $\widetilde L_2$ w.r.t. each other  and count the differences of multiplicities a.e with respect to each part according to Lemma \ref{proposition orthogonal sum}; we skip the details. 
\end{proof}

\subsection{Boundary relations and the Titchmarsh--Kodaira formula}

Using the abstract setting of boundary relations leads to a unified approach to the spectral theory of many concrete operators. 
A recent standard reference for this theory is \cite{Behrndt-Hassi-deSnoo-2020}; we shall sometimes also refer to 
\cite{Derkach-Hassi-Malamud-deSnoo-2006}. Let us now recall some basic facts used in the present paper.

\begin{definition}
	Let $H$ and $B$ be Hilbert spaces, $S$ be a closed symmetric linear relation in $H$ and $\Gamma$ be a linear relation from $H^2$ to $B^2$ (i.e., $\Gamma\subseteq H^2\times B^2$). 
	Then $\Gamma$ is called a \emph{boundary relation} for $S^*$, if the following holds.
	\begin{enumerate}[(BR1)]
	\item The domain of $\Gamma$ is contained in $S^*$ and is dense there.
	\item For all $((f;g);(a;b)),((f';g');(a';b'))\in\Gamma$ the \emph{abstract Green's identity} holds:
		$$
		(g,f')_{H}-(f,g')_{H}=(b,a')_{B}-(a,b')_{B}
		.
		$$
	\item The relation $\Gamma$ is maximal with respect to the properties {\rm(BR1)} and {\rm(BR2)}.
	\end{enumerate}
	If $\mul\Gamma\cap\big(\{0\}\times B\big)=\{0\}$, then we say that $\Gamma$ is of \emph{function type}. 
	If $\mul\Gamma=\{(0;0)\}$, then $\Gamma$ is called a \emph{boundary function}.
\end{definition}

In what follows we consider only boundary relations with $B=\mathbb C^n$. In this case it is known that deficiency indices of
$S$ are equal and the domain of $\Gamma$ is equal to $S^*$. If additionally $\mul\Gamma=\{(0;0)\}$, then $\Gamma$ is a bounded
operator from $S^*$ to $\mathbb C^{2n}$ which can be written in the form $\Gamma=(\Gamma_1,\Gamma_2)$, $\Gamma_i:S^*\to\bb C^n$,
$i=1,2$, and then the collection $(\bb C^n;\Gamma_1,\Gamma_2)$ is called an {\it (ordinary) boundary triplet} for $S^*$. 

For the case of a boundary triplet $(\bb C^n;\Gamma_1,\Gamma_2)$ for $S^*$ symmetric and selfadjoint extensions of $S$ can be described in a neat way. 
To this end introduce the indefinite scalar product $[\cdot,\cdot]$ on $\bb C^n\times\bb C^n$ defined by the Gram operator 
\[
	J_{\bb C^n}:=i
	\begin{pmatrix}
		0 & I_{\mathbb C^n}
		\\
		-I_{\mathbb C^n} & 0
	\end{pmatrix}
	.
\]
Explicitly, this is 
\[
	\big[(a;b),(a';b')\big]:=i\big((b,a')_{\bb C^n}-(a,b')_{\bb C^n}\big)
	.
\]
Then symmetric extensions $A\subseteq\mc H\times\mc H$ of $S$ correspond bijectively  to $[\cdot,\cdot]$--neutral subspaces 
$\theta\subseteq\bb C^n\times\bb C^n$ by means of the inverse image map $\Gamma^{-1}$. In this correspondence $A$ is selfadjoint
if and only if $\theta$ is a maximal neutral subspace (equivalently, $\theta$ is neutral and $\dim\theta=n$). 
Such subspaces are sometimes also called \emph{Lagrange planes}, e.g., \cite{Harmer-2000}. One has $S=\ker\Gamma$ and the
induced linear operator $\widetilde \Gamma$ on the quotient linear space $S^*/S$ is its linear isomorphism with $\bb C^{2n}$. In
the general case of a boundary relation $\Gamma$ there is a linear isomorphism between the quotient spaces $S^*/S={\rm \dom\,}\Gamma/\ker\Gamma$ and $\ran\Gamma/{\rm \mul\,}\Gamma$.

\begin{definition}
	Let $\Gamma$ be a boundary relation. The map $M$ which assigns to a point $z\in\bb C\setminus\bb R$ the linear relation
	$$
		M(z):=\big\{(a;b)\in\bb C^n\times\bb C^n:\ \exists\,f\in H\text{ with }\big((f;zf);(a;b)\big)\in\Gamma\big\}
	$$
	is called the \emph{Weyl family} of $\Gamma$.	
\end{definition}

\noindent
The Weyl family $M$ of a boundary relation of function type (where $B=\mathbb C^n$) is a $n\times n$-matrix-valued Herglotz 
function, and is also called the {\it Weyl function}.

The Weyl family of a boundary relation is intimitely related to the spectral theory of selfadjoint extensions of $S$,
specifically, to the spectrum of the selfadjoint (see \cite[Proposition 4.22]{Derkach-Hassi-Malamud-deSnoo-2006}) relation
\begin{equation}
\label{a1}
	L:=\ker\big[\pi_1\circ\Gamma\big]
	,
\end{equation}
where $\pi_1:\bb C^n\times\bb C^n\to\bb C^n$ is the projection onto the first component $\pi_1((a;b)):=a$. In the case of a boundary triplet $(\bb C^n;\Gamma_1,\Gamma_2)$ we have $L=\ker\Gamma_1$.
Since the Weyl family comprises the information given by defect elements, naturally, a selfadjoint part of $S$ (including its multivalued part) cannot be 
accessed using $M(z)$. For this reason a statement about spectrum can only be expected for boundary relations of 
simple symmetric linear relations (which therefore are operators, but may be nondensely defined, so that their adjoint may be a multivalued linear relation). A cornerstone in Weyl theory is the \emph{Titchmarsh--Kodaira formula}.
We give a formulation which is in fact a generalization of the Titchmarsh--Weyl--Kodaira theory 
\cite{Titchmarsh-1962}, \cite{Weyl-1910}, \cite{Kodaira-1949} for one-dimensional Schr\"odinger operators to the abstract setting 
of boundary relations.

\begin{theorem}[Titchmarsh--Kodaira formula]
\label{a3}
	Let $S$ be a closed simple  symmetric linear relation in a Hilbert space $H$ and 
	$\Gamma\subseteq H^2\times\mathbb C^{2n}$ be a boundary relation of function type for $S^*$. Consider the selfadjoint extension $L:=\ker[\pi_1\circ\Gamma]$ of
	$S$. Moreover, let $M$ be the Weyl function associated with $\Gamma$, let 
	$\Omega$ be the measure in the Herglotz integral representation of $M$, and let $\rho$ be its trace measure 
	$\rho:=\tr\Omega$. 

	Then the operator part of $L$ is unitarily equivalent to the operator of multiplication by independent variable in the space 
	$L_2(\bb R,\Omega)$. The spectral multiplicity function $N_L$ of $L$ is given as 
	\begin{equation}\label{p3}
		N_L(x)=\rank\frac{d\Omega}{d\rho}(x)\text{ for }\rho\text{-a.a. }x\in\bb R.
	\end{equation}
\end{theorem}

\subsection{Pasting of boundary relations}\label{subsection pasting}

In this subsection we describe what we understand by a pasting of boundary relations by means of interface conditions. 
For more details on operations with boundary relations we refer the reader to \cite{Behrndt-Hassi-deSnoo-2020} or 
\cite[\S3]{Simonov-Woracek-2014}. We restrict ourselves to pastings of relations from a particular simple class ($B=\mathbb C$),
however the construction that we use would also make sense in a more general case.

We use the word {\it pasting} in two meanings: for boundary relations $\Gamma_l$ and for selfadjoint linear relations $L_l=\ker[\pi_1\circ\Gamma_l]$. For the former a pasting is obtained from a fractional linear transform $w$ of $\Gamma_0=\prod_{l=1}^n\Gamma_l$, where the matrix $w$ is $J$-unitary and is constructed in a non-unique way from matrices $A$ and $B$ which determine the interface condition. For the latter the pasting is uniquely determined by $A$ and $B$.

We consider the following setting.
	\begin{enumerate}[(D1)]
		\item Let $n\geq 2$ and for $l\in\{1,\ldots,n\}$ let either $H_l$ be a Hilbert space, $S_l$ a simple closed
		symmetric linear relation in $H_l$ with deficiency index $(1,1)$ (which is hence an operator, possibly not
		densely defined), $\Gamma_l\subseteq H_l^2\times\bb C^2$ a boundary relation of function type for $S_l^*$ and $L_l=\ker[\pi_1\circ\Gamma_l]$ a selfadjoint 
		linear relation, or (an \emph{artificial edge}) $H_l=\{0\}$, $S_l=L_l=S_l^*=\{(0;0)\}$, and $\Gamma_l=\{0\}^2\times\mul\Gamma_l\subseteq\{0\}^2\times\bb C^2$.
		We assume that for at least one $l$ we do not have an artificial edge. 
		\item 
		Denote $H:=\prod_{l=1}^nH_l$, $S:=\prod_{l=1}^nS_l$, $\Gamma_0:=\prod_{l=1}^n\Gamma_l$ and $L_0:=\prod_{l=1}^nL_l$.
		\item
		Let $A,B\in\bb C^{n\times n}$ be such that $AB^*=BA^*$ and $\rank(A,B)=n$.
	\end{enumerate}

Note that our assumptions in {\rm(D1)} imply that boundary relations $\Gamma_l$ which are not boundary functions must be
artificial. Namely, if $\mul\Gamma_l\neq\{0\}^2$, then by a dimension argument we must have $H_l=\{0\}$ and
$\Gamma_l=\{0\}^2\times\{(a;m_la),a\in\mathbb C\}$ with some constant $m_l\in\bb R$ serving as Weyl function. Allowing such degenerate cases has 
meaningful applications: in \cite{Simonov-Woracek-2014} using a construction with an artificial edge we showed that
Aronszajn--Donoghue result \cite{Aronszajn-1957}, \cite{Donoghue-1965} can be deduced from the result of that work
(these results actually imply each other).

We see that $S$ is a simple closed symmetric linear relation in $H$ with deficiency index at most $(n,n)$, and
$$
\Gamma_0=\Big\{
	\big((f;g);(a;b)\big)\in H^2\!\times\!(\bb C^n\!\times\!\bb C^n):\ \big((f_l;g_l);(a_l;b_l)\big)\in\Gamma_l, l=1,\dots,n
	\Big\}
$$
is a boundary relation for $S^*$. The Weyl family of $\Gamma_0$ is given as the diagonal relation $M_0=\diag(m_1,\ldots,m_n)$,
where $m_l$ are Weyl functions of boundary relations $\Gamma_l$ or the real constants associated with artificial edges,
respectively. Since $\Gamma_l$ are of function type, so is $\Gamma_0$, and
$M_0$ is a diagonal matrix function. Obviously
\[
	L_0=\ker[\pi_1\circ\Gamma_0]
	,
\]
and we think of $L_0$ as the \emph{uncoupled} selfadjoint relation. 

In order to model interaction between the boundary relations $\Gamma_1,\ldots,\Gamma_n$, one can use fractional linear transforms defined by $J_{\bb C^n}$-unitary matrices $w$. Recall here that a matrix $w\in\bb C^{2n\times 2n}$ is called $J_{\bb C^n}$-unitary, if $w^*J_{\bb C^n}w=J_{\bb C^n}$. Also note that $w$ is $J_{\bb C^n}$-unitary if and only if $w^*$ is:
\begin{multline*}
	w^*J_{\bb C^n}w=J_{\bb C^n}\ \Leftrightarrow\ w^*(iJ_{\bb C^n})\cdot w(iJ_{\bb C^n})=I_{\bb C^n}
	\\
	\ \Leftrightarrow\ w(iJ_{\bb C^n})\cdot w^*(iJ_{\bb C^n})=I_{\bb C^n}\ \Leftrightarrow\ wJ_{\bb C^n}w^*=J_{\bb C^n}
	.
\end{multline*}
Given matrices $A,B$ which satisfy (D3) we can construct a $J_{\bb C^n}$-unitary matrix $w\in\bb C^{2n\times 2n}$ with $A$ and $B$ used as upper blocks. The following result is of a folklore kind, however, for completeness we provide its proof.

\begin{lemma}\label{lemma J-unitarity}
	Let $A,B\in\bb C^{n\times n}$ be given. There exist $C,D\in\bb C^{n\times n}$ such that the matrix 
	\begin{equation}\label{p4}
	w:=\begin{pmatrix}A&B\\C&D\end{pmatrix}
	\end{equation}
	is $J_{\bb C^n}$-unitary, if and only if $AB^*=BA^*$ and $\rank(A,B)=n$. 
\end{lemma}
\begin{proof}
	Multiplying out the product $wJ_{\bb C^n}w^*$, shows that $w$ is $J_{\bb C^n}$-unitary if and only if the four 
	equations 
	\begin{equation}
	\label{a2}
	BA^*-AB^*=0,\ DA^*-CB^*=I_{\bb C^n},\ BC^*-AD^*=-I_{\bb C^n},\ DC^*-CD^*=0
	,
	\end{equation}
	hold. 
	
	The backwards implication readily follows: if we find $C$ and $D$ such that \eqref{p4} is $J_{\bb C^n}$-unitary, then $BA^*=AB^*$ and $\ker[(A,B)^*]=\{0\}$. In order to show the forward implication, assume that $A$ and $B$ satisfy these two conditions. Then the matrix 
	\[
	X:=AA^*+BB^*=(A, B) \begin{pmatrix} A^*\\ B^* \end{pmatrix}
	\]
	is positive definite. Set $C:=-X^{-1}B$ and $D:=X^{-1}A$. Then all four relations in \eqref{a2} are fulfilled.
\end{proof}

\noindent
Given $w$ of the form \eqref{p4}, consider the relation
\[
	\Gamma_w:=w\circ\Gamma_0=
	\Big\{\big((f;g);(Aa+Bb;Ca+Db)\big):\big((f;g);(a;b)\big)\in\Gamma_0\Big\}
	.
\]
One can show that $\Gamma_w$ is a boundary relation for $S^*$ and that the Weyl family of $\Gamma_w$ is
\begin{align*}
	M_w(z)= &\, \Big\{(Aa+Bb;Ca+Db):(a;b)\in M_0(z)\Big\}
	\\
	= &\, \Big\{\big((A+BM_0(z))a;(C+DM_0(z))a\big):a\in\bb C^n\Big\},\quad z\in\bb C\setminus\bb R
	.
\end{align*}
If $A+BM_0(z)$ is invertible for every $z\in\bb C\setminus\bb R$, then 
\begin{equation}\label{p7}
	M_w(z)=(C+DM_0(z))(A+BM_0(z))^{-1},\quad z\in\bb C\setminus\bb R
	,
\end{equation}
is a matrix valued function and $\Gamma_w$ is of function type.

The selfadjoint relation \eqref{a1} for $\Gamma_w$ is given as 
\begin{multline}\label{L_AB}
	L_{A,B}:=\ker(\pi_1\circ\Gamma_w)
	\\
	=\Big\{(f;g)\in H^2:\ \exists (a;b)\in\bb C^{2n}:
	\big((f;g);(a;b)\big)\in\Gamma_0,\ Aa+Bb=0\Big\}
	.
\end{multline}
It depends only on the first row of $w$, i.e., on the matrices $A$ and $B$,  and hence we may legitimately call $L_{A,B}$ the
\emph{pasting} of the relations $L_l=\ker(\pi_1\circ\Gamma_l)$ by means of the interface conditions $(A,B)$. We think of
$L_{A,B}$ as the \emph{coupled} selfadjoint relation. 

The relation $L_{A,B}$ is a finite-rank perturbation of $L_0$ in the resolvent sense with the actual rank of the perturbation
not exceeding $n$. This holds simply because both are extensions of the symmetry $S$ which has deficiency index at most
$(n,n)$. The following lemma helps to estimate this rank.

Denote for $A,B\in\bb C^{n\times n}$
\begin{equation}\label{theta}
	\theta_{A,B}:=\{(a,b)\in\bb C^n\times\bb C^n:\ Aa+Bb=0\}.
\end{equation}
\begin{lemma}\label{lemma ranks}
	Let $A,B\in\bb C^{n\times n}$ and $A',B'\in\bb C^{n\times n}$ satisfy assumption {\rm (D3)}. Then
	$$
	\rank\big[(L_{A,B}-zI)^{-1}-(L_{A',B'}-zI)^{-1}\big]\leq\dim\Big(\raisebox{3pt}{$\theta_{A,B}$}\Big/\raisebox{-3pt}{$\theta_{A,B}\cap\theta_{A',B'}$}\Big),\quad z\in\bb C\setminus\bb R
	.
	$$
\end{lemma}	

\begin{proof}
	Under (D3), $\theta_{A,B}$ and $\theta_{A',B'}$ are Lagrange planes in $\bb C^n\times\bb C^n$ which correspond to the relations $L_{A,B}$ and $L_{A',B'}$ in the sense that $\Gamma_0^{-1}(\theta_{A,B})=L_{A,B}$ and $\Gamma_0^{-1}(\theta_{A',B'})=L_{A',B'}$. For each $z\in\bb C\setminus\bb R$ the rank of the resolvent difference is 
	$$
	\rank\big[(L_{A,B}-zI)^{-1}-(L_{A',B'}-zI)^{-1}\big]
	=\dim\Big(\raisebox{3pt}{$L_{A,B}$}\Big/\raisebox{-3pt}{$L_{A,B}\cap L_{A',B'}$}\Big)
	.
	$$
	We have:
	%\Big(\raisebox{3pt}{$$}\Big/\raisebox{-3pt}{$$}\Big)
	\begin{multline*}
	\dim\Big(\raisebox{3pt}{$L_{A,B}$}\Big/\raisebox{-3pt}{$L_{A,B}\cap L_{A',B'}$}\Big)
	\\
	=\dim\Big(\raisebox{3pt}{$L_{A,B}$}\Big/\raisebox{-3pt}{$S$}\Big)
	-\dim\Big(\raisebox{3pt}{$L_{A,B}\cap L_{A',B'}$}\Big/\raisebox{-3pt}{$S$}\Big)
	\\
	\leq\dim\Big(\raisebox{3pt}{$\Gamma_0^{-1}(\theta_{A,B})$}\Big/\raisebox{-3pt}{$\ker\Gamma_0$}\Big)
	-\dim\Big(\raisebox{3pt}{$\Gamma_0^{-1}(\theta_{A,B}\cap\theta_{A',B'})$}\Big/\raisebox{-3pt}{$\ker\Gamma_0$}\Big)
	\\
	=\dim\Big(\raisebox{3pt}{$\theta_{A,B}\cap\ran\Gamma_0$}\Big/\raisebox{-3pt}{$\theta_{A,B}\cap\mul\Gamma_0$}\Big)
	\\
	-\dim\Big(\raisebox{3pt}{$\theta_{A,B}\cap\theta_{A',B'}\cap\ran\Gamma_0$}\Big/\raisebox{-3pt}{$\theta_{A,B}\cap\theta_{A',B'}\cap\mul\Gamma_0$}\Big)
	\\
	\leq
	\dim\Big(\raisebox{3pt}{$\theta_{A,B}\cap\ran\Gamma_0$}\Big/\raisebox{-3pt}{$\theta_{A,B}\cap\theta_{A',B'}\cap\mul\Gamma_0$}\Big)
	\\
	-\dim\Big(\raisebox{3pt}{$\theta_{A,B}\cap\theta_{A',B'}\cap\ran\Gamma_0$}\Big/\raisebox{-3pt}{$\theta_{A,B}\cap\theta_{A',B'}\cap\mul\Gamma_0$}\Big)
	\\
	=\dim\Big(\raisebox{3pt}{$\theta_{A,B}\cap\ran\Gamma_0$}\Big/\raisebox{-3pt}{$\theta_{A,B}\cap\theta_{A',B'}\cap\ran\Gamma_0$}\Big)
	\\
	\leq
	\dim\Big(\raisebox{3pt}{$\theta_{A,B}$}\Big/\raisebox{-3pt}{$\theta_{A,B}\cap\theta_{A',B'}$}\Big).
	\end{multline*}
\end{proof}
	
\begin{corollary}\label{lemma rank B}
	Let $A,B\in\bb C^{n\times n}$ satisfy {\rm (D3)}. Then
	$$
	\rank\big[(L_{A,B}-zI)^{-1}-(L_0-zI)^{-1}\big]\leq\rank B,\quad z\in\bb C\setminus\bb R
	.
	$$
\end{corollary}	
	
\begin{proof}
	We have	$\theta_0:=\theta_{I,0}=\{(a,b)\in\bb C^n\times\bb C^n:\ a=0\}$. By the lemma for $z\in\bb C\setminus\bb R$
	$$
		\rank\big[(L_{A,B}-zI)^{-1}-(L_{0}-zI)^{-1}\big]\leq
		\dim\Big(\raisebox{3pt}{$\theta_{A,B}$}\Big/\raisebox{-3pt}{$\theta_{A,B}\cap\theta_0$}\Big)
		,
	$$
	and this dimension is computed as 
	$$
		\dim\Big(\raisebox{3pt}{$\theta_{A,B}$}\Big/\raisebox{-3pt}{$\theta_{A,B}\cap\theta_0$}\Big)
		=n-\dim(\theta_{A,B}\cap\theta_0)=n-\dim\ker B=\rank B.
	$$
\end{proof}

\section{Discussion of interface conditions}\label{section interface condition}

In our present work we investigate spectral properties of coupled operators $L_{A,B}$ whose interface conditions $(A,B)$ are
subject to an additional mixing condition. Namely, besides the general condition
\begin{enumerate}[(D3)]
\item $AB^*=BA^*$ and $\rank(A,B)=n$
\end{enumerate}
we are going to assume the condition
\begin{enumerate}[(D4)]
\stepcounter{enumi}
\item each set of $\rank B$ many different columns of $B$ is linearly independent.
\end{enumerate}
Note that this is a condition on $B$ only.

We have no fully precise intuition for {\rm(D4)}. However, in some sense it is related to how different edges are mixed. 
To illustrate the situation, observe that the matrix $B_{st}$ from the standard interface 
conditions in \eqref{standard ic} satisfies {\rm(D4)}, and apparently the Kirchhoff condition 
$\sum_{l=1}^n u_l'(0)=0$ for a star-graph combines all edges. 
On the other hand consider for instance the matrices
\[
	A_1:=
	\begin{pmatrix}
		1 & -1 & 0
		\\
		0 & 0 & 0
		\\
		0 & 0 & 1
	\end{pmatrix}
	,\quad
	A_2=
	\begin{pmatrix}
		1 & -1 & 0
		\\
		0 & 0 & 1
		\\
		0 & 1 & 0
	\end{pmatrix}
	,\quad
	B:=
	\begin{pmatrix}
		0 & 0 & 0
		\\
		1 & 1 & 0
		\\
		0 & 0 & 1
	\end{pmatrix}
	.
\]
Obviously, $B$ does not satisfy {\rm(D4)}. Thinking of a situation as in \eqref{i3}--\eqref{i4}, we see that interface conditions with 
this matrix $B$ fail to mix all edges in the second component of their boundary values. For example, the operator corresponding 
to the interface conditions $(A_1,B)$ splits in two uncoupled parts, one corresponding to the first two edges and another to the 
third edge. At the same time, the operator corresponding to $(A_2,B)$ will mix all edges.

The condition {\rm(D4)} can be characterized in different ways. 
Here we call a linear subspace $\mc L$ of $\bb C^n$ a coordinate plane, if $\mc L=\spn\{e_{i_1},\ldots,e_{i_l}\}$ with some
$1\leq i_1<i_2<\ldots<i_l\leq n$, and where $e_i$ denotes the $i$-th canonical basis vector in $\bb C^n$. 

\begin{lemma}\label{lemma 1}
	For a matrix $B\in\bb C^{n\times n}$ the following statements are equivalent. 
	\begin{enumerate}[$(1)$]
	\item $B$ satisfies {\rm(D4)}.
	\item For every coordinate plane $\mc L$ in $\mathbb C^n$ with $\dim\mc L\leq\rank B$, the restriction of $B$ to $\mc L$ is injective. 
	\item For every positive semidefinite diagonal matrix $X$ it holds that 
		\begin{equation}\label{equ 2}
			\rank(BXB^*)=\min\{\rank X,\rank B\}
			.
		\end{equation}
	\end{enumerate}
\end{lemma}
\begin{proof}
	Assume (1), i.e., that $B$ satisfies {\rm(D4)}, and let $\mc L=\spn\{e_{i_1},\ldots,e_{i_l}\}$ be a 
	coordinate plane with dimension $l\leq\rank B$. The range of $B|_{\mc L}$ is the linear span of the corresponding 
	columns of $B$, and hence has dimension $l$. Thus $B|_{\mc L}$ is injective, and we have (2). 

	Assume $(2)$ and let $X$ be a positive semidefinite diagonal matrix. Since $X\geq 0$, we have 
	\[
		\ker(BXB^*)=\ker(X^{\frac 12}B^*)
		,
	\]
	and hence 
	\begin{equation}\label{equ 3}
		\rank(BXB^*)=\rank(X^{\frac 12}B^*)=\rank(BX^{\frac 12})
		.
	\end{equation}
	Denote $r:=\rank B$ and $r':=\rank X$. Clearly, $\rank(BX^{\frac 12})\leq\min\{r,r'\}$. 
	Let $1\leq i_1<\ldots<i_{r'}\leq n$ be those indices for which the corresponding diagonal entry of $X$ is nonzero. 
	Then $\ran X^{\frac 12}=\spn\{e_{i_1},\ldots,e_{i_{r'}}\}=:\mc L$. 
	If $r'\leq r$, then $B$ acts injectively on $\mc L$, and hence $\rank(BX^{\frac 12})=r'$. 
	If $r'>r$, then $B$ acts injectively on $\spn\{e_{i_1},\ldots,e_{i_r}\}$, and hence $\rank(BX^{\frac 12})=r$.
	Put together, we arrive at the asserted formula \eqref{equ 2}, i.e., we have (3). 

	Finally, assume $(3)$ and let $1\leq i_1<\ldots<i_r\leq n$. Let $X$ be the diagonal matrix having diagonal
	entry $1$ in the $i_l$-th columns, $l=1,\ldots,r$, and $0$ otherwise. Then the linear span $\mc K$ 
	of the $i_1,\ldots,i_r$-th columns of $B$ is nothing but the range of $BX^{\frac 12}$. Remembering \eqref{equ 3}, we 
	obtain from \eqref{equ 2} that 
	\[
		\dim\mc K=\rank(BXB^*)=r
		.
	\]
	This means that the $i_1,\ldots,i_r$-th columns of $B$ are linearly independent and we have (1).
	
	Thus items (1), (2) and (3) are equivalent. 
\end{proof}

\noindent
The selfadjoint relation $L_{A,B}$ defined as pasting by means of the interface conditions $(A,B)$, cf.\ \eqref{L_AB}, does not
uniquely correspond to $(A,B)$. Clearly, if $Q\in\bb C^{n\times n}$ is invertible, then $L_{A,B}=L_{QA,QB}$. Another way
to modify $A$ and $B$ without essentially changing the corresponding operator is to simultaneously permute columns
of $A$ and $B$; this corresponds to ``renumerating the edges'': Let $\pi$ be a permutation of $\{1,\ldots,n\}$, and let $P$ be
the corresponding permutation matrix. Then the operator $L_{A,B}$ defined as pasting of 
$\Gamma_1,\ldots,\Gamma_n$ with interface conditions $(A,B)$ is unitarily equivalent to the operator built
from $\Gamma_{\pi(1)},\ldots,\Gamma_{\pi(n)}$ with $(AP,BP)$, and hence these two operators share all their spectral properties. 

We are going to use the above two operations to reduce interface conditions to a suitable normal form. 
For matrices $A,B,\tilde A,\tilde B\in\bb C^{n\times n}$ let us write $(A,B)\sim(\tilde A,\tilde B)$, if
there exist $Q,P\in\bb C^{n\times n}$ with $Q$ invertible and $P$ a permutation matrix, such that 
\[
	(\tilde A,\tilde B)=Q(A,B)\begin{pmatrix} P & 0 \\ 0 & P\end{pmatrix}
	.
\]
Clearly, $\sim$ is an equivalence relation. Two equivalent pairs $(A,B)$ and $(\tilde A,\tilde B)$ together do or do
not satisfy {\rm(D3)}. Further, for any matrices $B,Q,P\in\bb C^{n\times n}$ with $Q$ invertible and $P$ a permutation matrix, 
the matrices $B$ and $QBP$ together do or do not satisfy {\rm(D4)}. 

\begin{lemma}\label{lemma 2}
	Let $A,B\in\bb C^{n\times n}$ and assume that {\rm(D3)} and {\rm(D4)} hold. Denote $r:=\rank B$. Then there exist 
	$A_1\in\bb C^{(n-r)\times r}$ and $A_2\in\bb C^{r\times r}$ with $A_2=A_2^*$, such that 
	($I_{\bb C^k}$ denotes the identity matrix of size $k\times k$, and block matrices are understood with respect to the
	decomposition $\bb C^n=\bb C^{n-r}\times\bb C^r$)
	\[
		(A,B)\sim\left(
		\begin{pmatrix} I_{\bb C^{n-r}} & A_1 \\ 0 & A_2 \end{pmatrix}
		\,,\,
		\begin{pmatrix} 0 & 0 \\ -A_1^* & I_{\bb C^r} \end{pmatrix}
		\right)
		.
	\]
\end{lemma}

\begin{proof}
	By the definition of $r$ we find an invertible $Q_1\in\bb C^{n\times n}$ such that 
	\[
		Q_1\cdot B=\begin{pmatrix} 0 & 0 \\ B_{21} & B_{22} \end{pmatrix}
	\]
	with some blocks $B_{21}\in\bb C^{r\times(n-r)}$ and $B_{22}\in\bb C^{r\times r}$. 
	Since $\rank(Q_1A,Q_1B)=\rank(A,B)=n$, the first $n-r$ rows of $Q_1A$ are linearly independent. 
	Hence, we find an invertible $Q_2\in\bb C^{(n-r)\times(n-r)}$ and a permutation matrix $P\in\bb C^{n\times n}$
	such that 
	\[
		\begin{pmatrix} Q_2 & 0 \\ 0 & I_{\bb C^r} \end{pmatrix}\cdot Q_1A\cdot P=
		\begin{pmatrix} I_{\bb C^{n-r}} & A_{12} \\ 0 & A_{22} \end{pmatrix}
	\]
	with some blocks $A_{12}\in\bb C^{(n-r)\times r}$ and $A_{22}\in\bb C^{r\times r}$. 
	By {\rm(D4)}, the last $r$ columns of $Q_1BP$ are linearly independent. Equivalently, the right lower $r\times r$-block 
	$B_{22}'$ of $Q_1BP$ is invertible. Setting $Q_3:=(B_{22}')^{-1}$, we obtain 
	\[
		\begin{pmatrix} I_{\bb C^{n-r}} & 0 \\ 0 & Q_3 \end{pmatrix}\cdot Q_1B
		\cdot P
		=
		\begin{pmatrix} 0 & 0 \\ B_{21}' & I_{\bb C^r} \end{pmatrix}
	\]
	with some block $B_{21}'\in\bb C^{r\times(n-r)}$. 

	Putting together, we have 
	\[
		\begin{pmatrix} Q_2 & 0 \\ 0 & Q_3 \end{pmatrix}Q_1\cdot(A,B)\cdot 
		\begin{pmatrix} P & 0 \\ 0 & P \end{pmatrix}
		=\left(
		\begin{pmatrix} I_{\bb C^{n-r}} & A_{12} \\ 0 & Q_3A_{22} \end{pmatrix}
		\,,\,
		\begin{pmatrix} 0 & 0 \\ B_{21}' & I_{\bb C^r} \end{pmatrix}
		\right)
	\]
	From $AB^*=BA^*$, we obtain that $Q_3A_{22}=(Q_3A_{22})^*$ and $B_{21}'=-A_{12}^*$. 
\end{proof}

\noindent
\textbf{Conclusion:} When investigating spectral properties of pastings of boundary relations with interface conditions given by
matrices $A$ and $B$ subject to {\rm(D3)} and {\rm(D4)}, we may restrict attention without loss of generality to pairs $(A,B)$ 
of the form 
\begin{equation}\label{equ 1}
	A=\begin{pmatrix} I_{\bb C^{n-r}} & A_1 \\ 0 & A_2 \end{pmatrix}
	\,,\quad
	B=\begin{pmatrix} 0 & 0 \\ -A_1^* & I_{\bb C^r} \end{pmatrix}
	,
\end{equation}
with some $A_1\in\bb C^{(n-r)\times r}$ and $A_2\in\bb C^{r\times r}$, $A_2=A_2^*$. 

\begin{remark}\label{remark 1}
	Two pairs of the form \eqref{equ 1} can be equivalent modulo $\sim$. For example,
	\[
		\left(\begin{pmatrix} I_{\bb C^{n-r}} & A_1 \\ 0 & A_2 \end{pmatrix}
		\,,\,
		\begin{pmatrix} 0 & 0 \\ -A_1^* & I_{\bb C^r} \end{pmatrix}\right)
		\ \sim\ 
		\left(\begin{pmatrix} I_{\bb C^{n-r}} & \tilde A_1 \\ 0 & \tilde A_2 \end{pmatrix}
		\,,\,
		\begin{pmatrix} 0 & 0 \\ -\tilde A_1^* & I_{\bb C^r} \end{pmatrix}\right)
	\]
	if there exist permutation matrices $P_1\in\bb C^{(n-r)\times(n-r)}$ and $P_2\in\bb C^{r\times r}$ such that 
	\[
		\tilde A_1=P_1^*A_1P_2,\quad \tilde A_2=P_2^*A_2P_2
		.
	\]
\end{remark}

\noindent
Validity of {\rm(D4)} for a matrix $B$ of form \eqref{equ 1} can also be characterized in different ways.

\begin{lemma}\label{lemma 3}
	Let $B_1\in\bb C^{r\times(n-r)}$, and set 
	\[
		B:=\begin{pmatrix} 0 & 0 \\ B_1 & I_{\bb C^r} \end{pmatrix}
		.
	\]
	Then the following statements are equivalent.
	\begin{enumerate}[$(1)$]
	\item $B$ satisfies {\rm(D4)}.
	\item All minors of size $r$ of the matrix $(B_1,I_{\bb C^r})$ are nonzero.
	\item All minors of size $n-r$ of the matrix $(I_{\bb C^{n-r}},-B_1^*)$ are nonzero.
	\item All minors of $B_1$ are nonzero. 
	\end{enumerate}
\end{lemma}
\begin{proof}
	Since the upper $n-r$ rows of $B$ contain only zeros, a set of columns of $B$ is linearly independent if and only if
	the set consisting of the same columns of $(B_1,I_{\bb C^r})$ is linearly independent. 
	This shows the equivalence of $(1)$ and $(2)$. 

	Let $k\leq\min\{n-r,r\}$, and let $1\leq j_1<\ldots<j_k\leq n-r$ and $n-r<i_1<\ldots<i_k\leq n$. 
	Denote by $d$ the minor of size $k$ of the matrix $B_1$ obtained by selecting the columns $j_1,\ldots,j_k$ and the rows 
	$i_1,\ldots,i_k$ of $B$. 
	Let $1\leq j'_1<\ldots<j'_{n-r-k}\leq n-r$ and $n-r<i'_1<\ldots<i'_{r-k}\leq n$ be the complementary index sets, i.e., 
	\[
		\{j_1,\ldots,j_k\}\cap\{j'_1,\ldots,j'_{n-r-k}\}=\{i_1,\ldots,i_k\}\cap\{i'_1,\ldots,i'_{r-k}\}=\emptyset
		.
	\]
	Then the minor $d$ is up to a sign equal to the minor of size $r$ of $(B_1,I_{\bb C^r})$ made from the 
	$j_1,\ldots,j_k,n-r+i'_1,\ldots,n-r+i'_{r-k}$-th columns of this matrix. Further, it is up to a sign and complex conjugation equal to the minor of size $n-r$ of $(I_{\bb C^{n-r}},-B_1^*)$ made from the $j'_1,\ldots,j'_{n-r-k},n-r+i_1,\ldots,n-r+i_k$-th columns 
	of this matrix. This shows the equivalence of $(2)$, $(3)$, and $(4)$. 
\end{proof}

\section{Formulation of the main result}\label{section formulation}

In this section we formulate the main result of the present paper, the below Theorem~\ref{theorem main}. 
It is preceded by a corresponding theorem about the point spectrum, Theorem~\ref{theorem point}, 
which we state independently for several reasons: behavior of multiplicity of
eigenvalues may serve as a simple and elementary model for the behavior of spectral multiplicity of singular spectrum, 
it allows an independent proof by linear algebra, and one assumption from Theorem \ref{theorem main}, (D5), can be dropped. 

The result about the point spectrum now reads as follows. 

\begin{theorem}\label{theorem point}
	Let data be given according to {\rm(D1)}--{\rm(D4)} above. Let $L_{A,B}$ be the selfadjoint relation constructed by pasting $L_1,\ldots,L_n$ by means of the interface conditions $(A,B)$, cf.\ \eqref{L_AB}. Let $r:=\rank B$. For $x\in\bb R$ denote by $N^p_{A,B}(x)$ the multiplicity of $x$ as an eigenvalue of $L_{A,B}$, and by $N^p_0(x)$ the 
	multiplicity of $x$ as an eigenvalue of $L_0$, i.e., 
	\begin{align*}
		N^p_{A,B}(x):= &\, \dim\ker(L_{A,B}-xI)
		,
		\\
		N_0^p(x):= &\, \dim\ker(L_0-xI)=\#\{l\in\{1,\ldots,n\}:x\in\sigma_p(L_l)\}
		.
	\end{align*}
	Then the following statements hold.
	\begin{enumerate}[\rm(P1)]
	\item If $N_0^p(x)\geq r$, then $N^p_{A,B}(x)=N_0^p(x)-r$.
	\item If $N_0^p(x)< r$, then $N^p_{A,B}(x)\leq r-N_0^p(x)$.
	\end{enumerate}
\end{theorem}

\begin{remark}
	 Recall that for any choice of a representative of the multiplicity function $N_0$ we have $N_0(x)=N_0^p(x)$ 
	if $N_0^p(x)>0$ (the same for $N_{A,B}$). It is possible that $N_0^p(x)=0$ while $N_0(x)>0$.
\end{remark}

The formulation of the result for singular spectrum is quite more elaborate compared to point spectrum, 
but only because of the measure theoretic nature of the involved quantities. 
The basic message is quite the same. Again there occurs a distinction into two cases, many layers of spectrum vs.\ few layers of 
spectrum. If locally at a point the uncoupled operator $L_0$ has many layers of spectrum compared to the rank of matrix $B$, 
then after performing the perturbation the multiplicity will have decreased by $\rank B$. In particular, if we had exactly 
$\rank B$ many layers, the point will not anymore contribute to the spectral measure. 
If $L_0$ has few layers of spectrum, then the multiplicity will not become too large after performing the perturbation. 
Depending on the ratio of $r$ and $N_0(x)$ it may increase or must decrease.

We use the simplified notation $\mathds{1}_{\{N_0=l\}}\cdot\nu$ for the part $\mathds{1}_{Y_l}\cdot\nu$ of a measure $\nu\ll\mu$ (and their sums such as $\mathds{1}_{\{N_0>r\}}\cdot\nu$), cf. the definition of $Y_l$ and $N_0$ in Section \ref{section multiplicity}. This definition does not depend on the choice of a representative of the equivalence class of sets $\{N_0=l\}$ owing to absolute continuity of $\nu$ w.r.t. $\mu$.

\begin{theorem}\label{theorem main}
	Let data be given according to {\rm(D1)}--{\rm(D4)} above. Let $w$ be a $J$-unitary matrix provided by 
	Lemma \ref{lemma J-unitarity}, cf. \eqref{p4}, $\Gamma_w=w\circ\Gamma_0$ be the pasting of 
	$\Gamma_1$,\ldots,$\Gamma_n$, and
	$L_{A,B}=\ker[\pi_1\circ\Gamma_w]$ be the pasting of $L_1$,\ldots,$L_n$ by means of the interface condition $(A,B)$, cf.
	\eqref{L_AB}. Let $m_l$ be the (scalar) Weyl functions for $\Gamma_l$, $M_0$ and $M_w$ be the matrix valued Weyl functions for
	$\Gamma_0$ and $\Gamma_w$, and $\mu_l$, $\Omega_0$, $\Omega_w$ be the measures in their Herglotz integral representations  (if the ``$l$-th edge'' is ``artificial'', we assume $\mu_l=0$).
	Let $\mu:=\sum_{l=1}^n\mu_l=\tr\Omega_0$ and $\rho:=\tr\Omega_w$ be scalar spectral measures of the linear selfadjoint
	relations $L_0$, $L_{A,B}$, and $N_0$, $N_{A,B}$ be their spectral multiplicity functions.

	Let $\mu=\mu_{ac}+\mu_s$ and $\rho=\rho_{ac}+\rho_s$ be Lebesgue decompositions of $\mu$ and $\rho$ w.r.t.\ the Lebesgue measure $\lambda$, 
	and $\rho_s=\rho_{s,ac}+\rho_{s,s}$ the Lebesgue decomposition of $\rho_s$ w.r.t.\ $\mu$. Let $r:=\rank B$.
	
	Assume in addition that we do not have too many artificial edges in the sense that 
	\begin{enumerate}[\rm(D5)]
	\item 		$\#\{l:\mul\Gamma_l=\{(0;0)\}\}\geq r$.
	\end{enumerate}

	Then the following statements hold. 
	\begin{enumerate}[{\rm(S1)}]
		\item $\rho_{ac}\sim\mu_{ac}$ and $N_{A,B}=N_0$ holds $\rho_{ac}$-a.e.
		\item $\mathds{1}_{\{N_0=r\}}\cdot\rho_{s,ac}=0$ or, equivalently, $\rho_{s,ac}\sim\mathds{1}_{\{N_0>r\}}\cdot\rho_{s,ac}+\mathds{1}_{\{0<N_0<r\}}\cdot\rho_{s,ac}$.
		\item $\mathds{1}_{\{N_0>r\}}\cdot\rho_{s,ac}\sim\mathds{1}_{\{N_0>r\}}\cdot\mu_s$ and 
			$N_{A,B}=N_0-r$ holds $\mathds{1}_{\{N_0>r\}}\cdot\rho_{s,ac}$-a.e.
		\item $N_{A,B}\leq r-N_0$ holds $\mathds{1}_{\{0<N_0<r\}}\cdot\rho_{s,ac}$-a.e.
		\item $N_{A,B}\leq r$ holds $\rho_{s,s}$-a.e.
	\end{enumerate}
\end{theorem}

\begin{remark}
	\phantom{}
	\begin{enumerate}
	\item 
	Item (S1) follows from \cite{Gesztesy-Tsekanovskii-2000} and is included here only for completeness. 
	Item (S5) follows from Lemma \ref{lemma unknown} and (S1) since $\rho=(\rho_{ac}+\rho_{s,ac})+\rho_{s,s}$ is the 
	Lebesgue decomposition of $\rho$ w.r.t. $\mu$.
	\item 		
	Note that spectral multiplicity functions $N_0$ and $N_{A,B}$ are not  defined uniquely and can be changed on $\mu$- and $\rho$-zero sets, respectively. However, these sets of non-uniqueness are $\rho_{s,ac}$-zero, because $\rho_{s,ac}\ll \rho$ and $\rho_{s,ac}\ll\mu$, thus we can compare $N_{A,B}$ and $N_0$ in (S3) and (S4).
	\item
	Items (S2) and (S3) correspond to the ``many layers case'', while item (S4) is the ``few layers case''.
	\end{enumerate}
\end{remark}

\noindent
The proof of Theorem~\ref{theorem point} is by linear algebra and is given in Section~\ref{section point}. 
The proof of Theorem~\ref{theorem main} is given in Sections~\ref{section many layers} and \ref{section few layers}. 
For the many layers case we proceed via the Titchmarsh--Kodaira formula and the boundary behavior of the Weyl function, 
while the few layers case is settled by a geometric reduction.

\section{The point spectrum}\label{section point}

In this section we prove Theorem~\ref{theorem point}. 
Throughout the section let data be given as in this theorem, where $A,B$ are assumed to be in block form as in \eqref{equ 1}. 
Moreover, fix $x\in\bb R$. 

As a first step we study the eigenspaces $\ker(S_l^*-xI)$ for each $l\in\{1,\ldots,n\}$ separately. 
Recall that $\dim\ker(S_l^*-xI)\leq 1$ and $\dim\mul\Gamma_l\leq 1$. If $\mul\Gamma_l\neq\{(0;0)\}$ then $H_l=\{0\}$ 
and in particular $\sigma_p(S_l^*)=\emptyset$. Since $\Gamma_l$ is of function type, $\mul\Gamma_l$ cannot contain $(0;1)$. 

We have $x\in\sigma_p(S_l^*)$ if and only if there exist $f_l\in H_l\setminus\{0\}$ and $a_l,b_l\in\bb C$ with 
\[
	\big((f_l;xf_l);(a_l;b_l)\big)\in\Gamma_l
	.
\]
Assume that $x\in\sigma_p(S_l^*)$. Then the boundary values $a_l,b_l$ are uniquely determined by the element $f_l$, 
and $f_l$ itself is unique up to a scalar multiple. Fix a choice of $f_l$. Since $S_l$ is simple, we have for the corresponding
boundary values $(a_l;b_l)\neq 0$. Moreover, $x\in\sigma_p(L_l)$ if and only if $a_l=0$. 

We fix the following notation throughout: 
\begin{enumerate}[$(1)$]
\item If $x\in\sigma_p(L_l)$, let $\tilde f_l\in H_l\setminus\{0\}$ be the unique element with 
	\[
		\big((\tilde f_l;x\tilde f_l);(0;1)\big)\in\Gamma_l
		.
	\]
\item If $x\in\sigma_p(S_l^*)\setminus\sigma_p(L_l)$, let $\tilde f_l\in H_l\setminus\{0\}$ and $m_l\in\bb C$ be the unique 
	elements with 
	\[
		\big((\tilde f_l;x\tilde f_l);(1;m_l)\big)\in\Gamma_l
		.
	\]
\item If $\mul\Gamma_l\neq\{(0;0)\}$, let $m_l\in\bb C$ be the unique element such that 
	\[
		((0;0);(1;m_l))\in\Gamma_l
		.
	\]
\end{enumerate}
Note that, by the abstract Green's identity, we have $m_l\in\bb R$  in this case. 
Moreover, we denote:
\begin{enumerate}[$(1)$]
\setcounter{enumi}{3}
\item If $x\notin\sigma_p(S_l^*)$, set $\tilde f_l:=0$ (including the case $\mul\Gamma_l\neq\{(0;0)\}$).
\item If $x\in\sigma_p(L_l)$, or $\mul\Gamma_l=\{(0;0)\}$ and $x\notin\sigma_p(S_l^*)$, set $m_l:=0$. 
\end{enumerate}
In the second step we identify $\ker(L_{A,B}-xI)$. Set 
\begin{align*}
	J_p &\, :=\big\{l\in\{1,\ldots,n\}\colon x\in\sigma_p(L_l)\big\},
	\\
	J_p^* &\, :=\big\{l\in\{1,\ldots,n\}\colon x\in\sigma_p(S_l^*)\big\},
	\\
	J_m &\, :=\big\{l\in\{1,\ldots,n\}\colon \mul\Gamma_l\neq\{(0;0)\}\big\},
\end{align*}
and let $\hat f_l\in H$ be the element whose $l$-th coordinate is $\tilde f_l$ and all other coordinates are $0$. 
For a subset $J$ of $\{1,\ldots,n\}$ let $D_J$ be the diagonal matrix whose diagonal entry in the $i$-th column is 
$1$ if $i\in J$ and $0$ otherwise. Moreover, let $M$ be the diagonal matrix with diagonal entries $m_1,\ldots,m_n$. 

\begin{lemma}\label{lemma 4}
	Set 
	\[
		\gamma:=AD_{(J_p^*\setminus J_p)\cup J_m}+B(D_{J_p}+MD_{(J_p^*\setminus J_p)\cup J_m})
	\]
	and let $\Xi:\bb C^n\to H$ be the map (we write $c=(c_l)_{l=1}^n$)
	\[
		\Xi(c):=\sum_{l=1}^nc_l\hat f_l
		.
	\]
	Then $\ker(L_{A,B}-xI)=\Xi(\ker\gamma)$, and in particular                                                                                                                                                                                                                                                                                                                                                                                                                                                                                                                                                                                                                                                                                                                                                                                                                                                                                                                                                                                                                                                                                                                                            
	\[
		\dim\ker(L_{A,B}-xI)=\dim\ker\gamma-\dim\big(\ker\gamma\cap\ker\Xi\big)
		.
	\]
\end{lemma}
\begin{proof}
	We have 
	\[
		\ker(S^*-xI)=\prod_{l=1}^n\ker(S_l^*-xI)=
		\left\{\sum_{l=1}^nc_l\hat f_l\colon c\in\bb C^n\right\}=\ran\Xi
		.
	\]
	In this representation of $\ker(S^*-xI)$ the constants $c_l$ for $l\notin J_p^*$ are irrelevant, since $\hat f_l=0$ for 
	$l\notin J_p^*$. This changes when we turn to $\ker(L_{A,B}-xI)\subseteq\ker(S^*-xI)$, since the values of $c_l$ for $l\in J_m$ influence boundary values. For every $\hat f=\sum_{l=1}^nc_l\hat f_l\in\ker(S^*-xI)$ there exists (possibly not unique) $(a;b)\in\bb C^{2n}$ such that
	\begin{equation}\label{ker}
			\Big((\hat f;x\hat f);(a;b)\Big)\in\Gamma_0
	\end{equation}
	and $\hat f\in\ker(L_{A,B}-xI)$, if and only if there exists $(a;b)\in\bb C^{2n}$ such that $Aa+Bb=0$ additionally to \eqref{ker}. According to (1)--(3) it should be that
	\[
		(a;b)=\sum_{l\in J_p}c_l(0;e_l)+\sum_{l\in J_p^*\setminus J_p}c_l(e_l;m_le_l)+\sum_{l\in J_m}c_l(e_l;m_le_l),
	\]
	which means that
\begin{equation*}
a=\sum_{l\in (J_p^*\setminus J_p)\cup J_m}c_le_l,
\quad
b=\sum_{l\in J_p}c_le_l+\sum_{l\in (J_p^*\setminus J_p)\cup J_m}c_lm_le_l.
\end{equation*}
	Then $Aa+Bb=0$ is equivalent to 
$$\big[AD_{(J_p^*\setminus J_p)\cup J_m}+B(D_{J_p}+MD_{(J_p^*\setminus J_p)\cup J_m})\big]c=0.$$
	We conclude that 
	\begin{multline*}
\ker(L_{A,B}-xI)
\\
=\left\{\sum_{l=1}^n c_l\hat f_l:c\in\bb C^n\wedge
\big[AD_{(J_p^*\setminus J_p)\cup J_m}+B(D_{J_p}+MD_{(J_p^*\setminus J_p)\cup J_m})\big]c=0\right\}
\\
=\left\{\sum_{l=1}^nc_l\hat f_l\colon c\in\ker\gamma\right\}=\Xi(\ker\gamma)
		.
	\end{multline*}
\end{proof}

\noindent
In the third step we compute or estimate, respectively, the dimensions of $\ker\gamma$ and $\ker\gamma\cap\ker\Xi$. 

\begin{lemma}\label{lemma 5}
	\phantom{}
	\begin{enumerate}[$(1)$]
	\item If $\#(J_p^*\setminus J_p)\cup J_m\leq n-r$, then 
		\begin{align*}
			& \dim\ker\gamma=\max\big\{\#J_p-r,0\big\}+n-\#J_p^*-\#J_m
			,
			\\
			& \dim\big(\ker\gamma\cap\ker\Xi\big)=n-\#J_p^*-\#J_m
			.
		\end{align*}
	\item If $\#(J_p^*\setminus J_p)\cup J_m>n-r$, then 
		\begin{align*}
			& \dim\ker\gamma\leq r-\#J_p
			,
			\\
			& \dim\big(\ker\gamma\cap\ker\Xi\big)\geq n-\#J_p^*-\#J_m
			.
		\end{align*}
	\end{enumerate}
\end{lemma}

\begin{proof}
	Consider the case that $\#(J_p^*\setminus J_p)\cup J_m\leq n-r$. We show that 
	\begin{equation}\label{equ 4}
		\ker\gamma=\big\{c\in\bb C^n\colon D_{J_p^*\setminus J_p}c=D_{J_m}c=BD_{J_p}c=0\big\}
		.
	\end{equation}
	The inclusion ``$\supseteq$'' holds since $m_l=0$ for $l\notin(J_p^*\setminus J_p)\cup J_m$. 
	Let $c\in\ker\gamma$. Since $D_{1,\ldots,n-r}B=0$, it follows that 
	\[
		D_{1,\ldots,n-r}AD_{(J_p^*\setminus J_p)\cup J_m}c=0
		.
	\]
	The left side is a linear combination of at most $n-r$ columns of the matrix $(I_{\bb C^{n-r}},A_1)$, and 
	Lemma~\ref{lemma 3} implies that $D_{(J_p^*\setminus J_p)\cup J_m}c=0$ and hence $Mc=0$. From this we obtain 
	\[
		BD_{J_p}c=0
		.
	\]
	 Thus the inclusion ``$\subseteq$'' also holds and the proof of \eqref{equ 4} is complete. 

	The sets $J_p^*\setminus J_p$, $J_m$, and $J_p$, are pairwise disjoint, and each $r$ columns of $B$ are linearly
	independent. Thus \eqref{equ 4} implies 
	\begin{multline*}
		\dim\ker\gamma
		=
		n-(\#J_p^*-\#J_p)-\#J_m-\min\{\#J_p,r\}
		\\
		=
		\max\big\{\#J_p-r,0\big\}+\big(n-\#J_p^*-\#J_m\big)
		.
	\end{multline*}
	We have 
	\[
		\ker\Xi=\{c\in\bb C^n\colon D_{J_p^*}c=0\}
		,
	\]
	and see that 
	\[
		\ker\gamma\cap\ker\Xi=\big\{c\in\bb C^n\colon D_{J_p^*}c=D_{J_m}c=0\big\}
		.
	\]
	Thus $\dim(\ker\gamma\cap\ker\Xi)=n-\#J_p^*-\#J_m$.

	Consider now the case that $\#(J_p^*\setminus J_p)\cup J_m>n-r$. We show that 
	\begin{equation}\label{equ 5}
		\dim\ran\gamma\geq(n-r)+\#J_p
		.
	\end{equation}
	Choose $n-r$ indices in $(J_p^*\setminus J_p)\cup J_m$, and let $\mc R$ be the linear span of the corresponding columns
	of $\gamma$. Denote by $\pi_+$ the projection in $\bb C^n$ onto the first $n-r$ coordinates 
	($\pi_+=D_{\{1,\ldots,n-r\}}$). The image $\pi_+(\mc R)$ is
	the span of $n-r$ columns of the matrix $(I_{\bb C^{n-r}},A_1)$, and hence has dimension $n-r$. It follows that 
	\[
		\dim\mc R=n-r,\quad \mc R\cap\ker\pi_+=\{0\}
		.
	\]
	Let $\mc R'$ be the linear span of the columns of $\gamma$ corresponding to indices in $J_p$. For those indices columns
	of $\gamma$ are actually columns of $B$, and it follows that 
	\[
		\dim\mc R'=\#J_p,\quad \mc R'\subseteq\ker\pi_+
		.
	\]
	Note here that $\#J_p<r$. Since $\mc R+\mc R'\subseteq\ran\gamma$, we obtain 
	\[
		\dim\ran\gamma\geq\dim\mc R+\dim\mc R'=(n-r)+\#J_p
		.
	\]
	From this, clearly, $\dim\ker\gamma\leq r-\#J_p$. 
	Next, we have $\ker D_{J_p^*\cup J_m}\subseteq\ker\gamma\cap\ker\Xi$, and hence 
	\[
		\dim\big(\ker\gamma\cap\ker\Xi\big)\geq n-\#J_p^*-\#J_m
		.
	\]
\end{proof}

\noindent
It is easy to deduce Theorem~\ref{theorem point} from the above lemma. 

\begin{proof}[Proof of Theorem~\ref{theorem point}]
	Assume that $N_0^p(x)\geq r$, and note that $N_0^p(x)=\#J_p$. 
	Then case $(1)$ in Lemma~\ref{lemma 5} takes place, and it follows that 
	\[
		N_{A,B}^p(x)=\dim\ker(L_{A,B}-xI)=\#J_p-r
		.
	\]
	Assume now that $N_0^p(x)<r$. If in Lemma~\ref{lemma 5} case $(1)$ takes place, then $N_{A,B}^p(x)=0$. If case $(2)$ takes
	place, we obtain the bound
	\[
		N_{A,B}^p(x)\leq r-\#J_p-n+\#J_p^*+\#J_m\leq r-\#J_p
		.
	\]
\end{proof}

\section{The many layers case}\label{section many layers}

In the present section we prove assertions {\rm(S2)} and {\rm(S3)} from Theorem~\ref{theorem main}. 
The argument rests on boundary limit formula for matrix valued Nevanlinna functions: provided that the boundary relation $\Gamma_w$ is of function type, we employ Titchmarsh--Kodaira formula and Lemma~\ref{lemma 20} to obtain that 
\begin{equation}\label{equ 12}
	N_{A,B}(x)=\rank\frac{d\Omega_w}{d\rho}(x)=\rank\,\lim_{z\downarrow x}\frac{\Im M_w(z)}{\Im\tr M_w(z)}
	\text{ for }\rho\text{-a.a. }x\in\bb R
	.
\end{equation}

Throughout this section let data and notation be as in Theorem~\ref{theorem main}. In addition assume that the matrices 
$A,B$ describing the interface condition are of the form 
\[
	A=\begin{pmatrix} I_{\bb C^{n-r}} & A_1 \\ 0 & A_2\end{pmatrix}
	,\quad
	B=\begin{pmatrix} 0 & 0 \\ -A_1^* & I_{\bb C^r}\end{pmatrix}
	,
\]
with some $A_1\in\bb C^{(n-r)\times r}$ and $A_2\in\bb C^{r\times r}$ such that all minors of $A_1$ are nonzero and $A_2=A_2^*$. 
Recall that, by the discussion in Section~\ref{section interface condition}, this assumption is no loss in generality. 

The algebraic core of the proof is the following auxiliary proposition.

\begin{proposition}\label{proposition 1}
	Let $r<n$, $B_1\in\bb C^{r\times(n-r)}$ be such that all minors of $B_1$ are nonzero, and let 
	$X_1\in\bb C^{(n-r)\times(n-r)}$ and $X_2\in\bb C^{r\times r}$ be two positive semidefinite diagonal matrices such that 
	$\rank X_1+\rank X_2\geq r$. Then 
	\begin{enumerate}[$(1)$]
	\item the matrix $B_1X_1B_1^*+X_2$ is invertible, and
	\item ${\displaystyle 
		\rank\big(X_1-X_1B_1^*(B_1X_1B_1^*+X_2)^{-1}B_1X_1\big)=\rank X_1+\rank X_2-r
		.
		}$
	\end{enumerate}
\end{proposition}
\begin{proof}
	To show $(1)$ set
	\[
		B:=\begin{pmatrix} 0 & 0 \\ B_1 & I_{\bb C^r} \end{pmatrix}
		,\quad
		X:=\begin{pmatrix} X_1 & 0 \\ 0 & X_2 \end{pmatrix}
		.
	\]
	Then 
	\[
		BXB^*=\begin{pmatrix} 0 & 0 \\ 0 & B_1X_1B_1^*+X_2 \end{pmatrix}
		.
	\]
	Clearly, $\rank X=\rank X_1+\rank X_2$. By Lemma~\ref{lemma 3} matrix $B$ satisfies {\rm(D4)}, and by 
	Lemma~\ref{lemma 1} we have 
	\[
		\rank\big(B_1X_1B_1^*+X_2\big)=\rank(BXB^*)=\min\{\rank X,\rank B\}=r
		.
	\]
	This means that the $r\times r$ matrix $B_1X_1B_1^*+X_2$ is invertible. 

	The assertion in $(2)$ requires a slightly more elaborate argument. Set 
	\[
		H:=B_1X_1B_1^*+X_2,\quad G:=I_{\bb C^{n-r}}-X_1B_1^*H^{-1}B_1
		.
	\]
	To avoid confusion with dimensions let us here make explicit that 
	\[
		X_1,G:\bb C^{n-r}\to\bb C^{n-r},\quad X_2,H,H^{-1}:\bb C^r\to\bb C^r
		,
	\]
	\[
		B_1:\bb C^{n-r}\to\bb C^r,\quad B_1^*:\bb C^r\to\bb C^{n-r}
		.
	\]
	The essence of the argument are the following four relations:
	\begin{align}
		\label{equ 10}
		& X_2\cdot H^{-1}B_1=B_1\cdot G,\quad G\cdot X_1B_1^*=X_1B_1^*H^{-1}\cdot X_2,
		\\
		\label{equ 11}
		& X_1B_1^*\cdot H^{-1}B_1=I_{\bb C^{n-r}}-G,\quad H^{-1}B_1\cdot X_1B_1^*=I_{\bb C^r}-H^{-1}X_2.
	\end{align}
	The equalities in \eqref{equ 10} follow by plugging $X_2=H-B_1X_1B_1^*$, the first relation in \eqref{equ 11} is just the definition 
	of $G$ and the second is $H^{-1}H=I_{\bb C^r}$. 

	Together \eqref{equ 10} and \eqref{equ 11} show that 
	$H^{-1}B_1|_{\ker G}$ and $X_1B_1^*|_{\ker X_2}$
	are mutually inverse bijections between $\ker G$ and $\ker X_2$. In particular, this implies that
	\[
		\dim\ker G=\dim\ker X_2
		.
	\]
	The assertion in $(2)$ is now easily deduced. 
	We have $\ker(GX_1)=X_1^{-1}(\ker G)$ (the preimage of $\ker G$). The definition of $G$ shows that 
	$\ker G\subseteq\ran X_1$, and hence 
	\[
		X_1|_{X_1^{-1}(\ker G)}:X_1^{-1}(\ker G)\to\ker G
	\]
	is surjective. Clearly, $\ker X_1\subseteq X_1^{-1}(\ker G)$, and hence 
	\[
		\ker\big(X_1|_{X_1^{-1}(\ker G)}\big)=\ker X_1
		.
	\]
	Together this implies that 
	\begin{multline*}
		\dim\ker(GX_1)=\dim X_1^{-1}(\ker G)
		\\
		=\dim\ker\big(X_1|_{X_1^{-1}(\ker G)}\big)+\dim\ker G=\dim\ker X_1+\dim\ker X_2
		,
	\end{multline*}
	and hence that 
	\begin{multline*}
		\rank(GX_1)=(n-r)-\dim\ker(GX_1)
		\\
		=(n-r)-\dim\ker X_1-\dim\ker X_2=\rank X_1+\rank X_2-r
		.
	\end{multline*}
\end{proof}	

\noindent
We return to the setting of the theorem. Denote for $z\in\bb C\setminus\bb R$ 
\[
M_0(z)=:\begin{pmatrix} M_{11}(z) & 0 \\ 0 & M_{22}(z) \end{pmatrix}
.
\]
Writing out the block form of $A+BM_0(z)$ gives 
\[
	A+BM_0(z)=\begin{pmatrix} I_{\bb C^{n-r}} & A_1 \\ -A_1^*M_{11}(z) & A_2+M_{22}(z) \end{pmatrix}
	.
\]
Denote also 
\begin{equation}\label{D}
D(z):=A_2+M_{22}(z)+A_1^*M_{11}(z)A_1.
\end{equation}
We see that $D(z)$ is the Schur complement of the upper-left block of $A+BM_0(z)$. 

\begin{lemma}\label{lemma 6}
	Let {\rm (D1)--(\rm D5)} hold with $r<n$. Then $\Gamma_w$ is of function type and for every $z\in\bb C\setminus\bb R$
	\begin{enumerate}[$(1)$]
	\item $D(z)$ is invertible,
	\item $A+BM_0(z)$ is invertible,
		\begin{equation}\label{ml1}
			(A+BM_0(z))^{-1}=
			\begin{pmatrix}
				I_{\bb C^{n-r}}-A_1D(z)^{-1}A_1^*M_{11}(z) & -A_1D(z)^{-1} 
				\\ 
				D(z)^{-1}A_1^*M_{11}(z) & D(z)^{-1}
			\end{pmatrix}
		\end{equation}
		and
		\begin{equation}\label{ml2}
			\Im M_w(z)=((A+BM_0(z))^{-1})^*\Im M_0(z)(A+BM_0(z))^{-1}
			.
		\end{equation}
	\end{enumerate}
\end{lemma}
\begin{proof}
	The imaginary part of $D(z)$ is 
	\[
		\Im D(z)=\Im M_{22}(z)+A_1^*\Im M_{11}(z)A_1
		.
	\]
	By assumption (D5) at least $r$ many of the functions $m_l$ are not equal to a real constant.
	Thus 
	\[
		\rank\Im M_{11}(z)+\rank\Im M_{22}(z)=\rank\Im M_0(z)\geq r
		,
	\]
	and Proposition~\ref{proposition 1}, $(1)$, implies that $\Im D$ is invertible. Since $\Im D(z)\geq 0$,  it follows that
	$D(z)$ is also invertible: if $D(z)u=0$, then 
	\[
		\|(\Im D)^{\frac 12}u\|^2=\big((\Im D)u,u\big)=\frac 1{2i}\big((Du,u)-(u,Du)\big)=0
		,
	\]
	and hence $u=0$. Since  the Schur complement $D(z)$ is invertible, 
	we conclude that $A+BM_0(z)$ is invertible and the formula \eqref{ml1} holds. 

	As we have noted in Subsection \ref{subsection pasting}, invertibility of $A+BM_0(z)$ implies that $M_w(z)$ is a matrix for every $z\in\bb C\setminus\bb R$ and that $\Gamma_w$ is of function type. Formula \eqref{ml2} holds by \cite[Lemma~6.3,~(i)]{Gesztesy-Tsekanovskii-2000}. 
\end{proof}

\begin{remark}
	The case $r=n$ needs separate attention (note that then all functions $m_l$ are not real constants). In this case there is no block structure and we can write
	$$
	A=A_{22},\quad B=I,\quad M_0(z)=M_{22}(z),\quad A+BM_0(z)=A_{22}+M_0(z)
	$$
	and
	$$
	\Im(A+BM_0(z))=\Im M_0(z)
	$$
	is invertible, hence $A+BM_0(z)$ is also invertible and $\Gamma_w$ is of function type, \eqref{ml1} becomes $(A+BM_0(z))^{-1}=(A_{22}+M_0(z))^{-1}(z)$ and \eqref{ml2} continues to hold.
\end{remark}

\noindent
For points $x\in\bb R$ with sufficiently nice properties formulae \eqref{ml1}--\eqref{ml2} can be used to compute \eqref{equ 12}. 
Here, and in the following, we denote $m:=\sum_{l=1}^n m_l$. 

\begin{lemma}\label{lemma 9}
	Let (D1)--(D5) hold with $r<n$ and let $x\in\bb R$. Assume that 
	\begin{enumerate}[$(1)$]
	\item the symmetric derivative $\frac{d\Omega_0}{d\mu}(x)$ exists,
	\item the symmetric derivative $\frac{d\mu}{d\lambda}(x)$ exists and is equal $\infty$, 
	\item $\lim\limits_{z\downarrow x}\frac{M_0(z)}{m(z)}$ exists and equals to $\frac{d\Omega_0}{d\mu}(x)$.
	\end{enumerate}
	 Set 
	$$
	\alpha(x):=\frac{d\Omega_0}{d\mu}(x),
	$$
	let as usual $\alpha_{11}(x)$ and $\alpha_{22}(x)$ denote 
	the diagonal blocks of $\alpha(x)$ of size $n-r$ and $r$, respectively,
	and let
	$$
	H(x):=\alpha_{22}(x)+A_1^*\alpha_{11}(x)A_1.
	$$
	If $\rank\alpha(x)> r$, then $H(x)$ is invertible, and the following limit exists, equals 
	\[
		\lim_{z\downarrow x}\frac{\Im M_w(z)}{\Im m(z)}=
		\begin{pmatrix}
			\alpha_{11}(x)\!-\!\alpha_{11}(x)A_1
			H(x)^{-1}A_1^*\alpha_{11}(x)
			& 0
			\\[2mm]
			0 & 0
		\end{pmatrix}
	\]	
	and its rank equals $\rank\alpha(x)-r$.
	Assume that in addition to $(1)$--$(3)$ it holds that 
	\begin{enumerate}[$(1)$]
	\setcounter{enumi}{3}
	\item the symmetric derivative $\frac{d\rho}{d\mu}(x)$ exists and is finite. 
	\end{enumerate}
	Then $\frac{d\rho}{d\mu}(x)=0$ if and only if $\rank\alpha(x)=r$; if $\rank\alpha(x)>r$, then the limit $\lim_{z\downarrow x}\frac{\Im M_w(z)}{\Im\tr M_w(z)}$ exists and
	\begin{equation}\label{equ 6}
		\rank\,\lim_{z\downarrow x}\frac{\Im M_w(z)}{\Im\tr M_w(z)}=\rank\alpha(x)-r
		.
	\end{equation}
\end{lemma}

\begin{proof}
	Let $(1)$--$(3)$ hold. Set $H(x):=\alpha_{22}(x)+A_1^*\alpha_{11}(x)A_1$. From (2) it follows by 
	Proposition \ref{proposition 20}, (2), that $m(z)\to\infty$ as $z\downarrow x$. From (3)  and \eqref{D} we have $\lim_{z\downarrow x}\frac{D(z)}{m(z)}=H(x)$.
	Proposition~\ref{proposition 1}, $(1)$, tells us that  if $\rank\alpha(x)>r$, then $H(x)$ is invertible. Hence,  in this case we may conclude that 
	\[
		\lim_{z\downarrow x}m(z)D(z)^{-1}=H(x)^{-1}
		.
	\]
	Now the representation \eqref{ml1} yields 
	\begin{align*}
		& \lim_{z\downarrow x}\big(A+BM_0(z)\big)^{-1}
		\\
		& =\lim_{z\downarrow x}
		\begin{pmatrix}
			I_{\bb C^{n-r}}-A_1\cdot m(z)D(z)^{-1}\cdot A_1^*\cdot 
			\frac{M_{11}(z)}{m(z)}
			&
			-A_1\cdot m(z)D(z)^{-1}\cdot\frac 1{m(z)}
			\\[2mm]
			m(z)D(z)^{-1}\cdot A_1^*\cdot\frac{M_{11}(z)}{m(z)}
			&
			m(z)D(z)^{-1}\cdot\frac 1{m(z)}
		\end{pmatrix}
		\\[2mm]
		& =
		\begin{pmatrix}
			I_{\bb C^{n-r}}-A_1H(x)^{-1}A_1^*\alpha_{11}(x) & 0
			\\[2mm]
			H(x)^{-1}A_1^*\alpha_{11}(x) & 0
		\end{pmatrix}
		.
	\end{align*}
	Since the symmetric derivative  $\frac{d\Omega_0}{d\mu}(x)=\alpha(x)$ exists and $\frac{d\mu}{d\lambda}(x)=\infty$, we  get by
	Proposition \ref{proposition 20}, (1), together with Lemma \ref{lemma matrix Herglotz} that
	\[
		\lim_{z\downarrow x}\frac{\Im M_0(z)}{\Im m(z)}=\alpha(x)
		.
	\]
	Using \eqref{ml2} we see that the limit exists
	\begin{align*}
		\lim_{z\downarrow x}\frac{\Im M_w(z)}{\Im m(z)}
		= &\, 
		\begin{pmatrix}
			I_{\bb C^{n-r}}-A_1H(x)^{-1}A_1^*\alpha_{11}(x) & 0
			\\[2mm]
			H(x)^{-1}A_1^*\alpha_{11}(x) & 0
		\end{pmatrix}
		^*
		\cdot
		\\[2mm]
		&\, \cdot
		\begin{pmatrix} \alpha_{11}(x) & 0 \\ 0 & \alpha_{22}(x) \end{pmatrix}
		\cdot
		\begin{pmatrix}
			I_{\bb C^{n-r}}-A_1H(x)^{-1}A_1^*\alpha_{11}(x) & 0
			\\[2mm]
			H(x)^{-1}A_1^*\alpha_{11}(x) & 0
		\end{pmatrix}
		\\[2mm]
		= &\, 
		\begin{pmatrix}
			\alpha_{11}(x)\!-\!\alpha_{11}(x)A_1
			\big(\alpha_{22}(x)\!+\!A_1^*\alpha_{11}(x)A_1\big)^{-1}A_1^*\alpha_{11}(x)
			& 0
			\\[2mm]
			0 & 0
		\end{pmatrix}
		.
	\end{align*}
 Applying Proposition~\ref{proposition 1}, $(2)$, we obtain that 
	\[
		\rank\,\lim_{z\downarrow x}\frac{\Im M_w(z)}{\Im m(z)}=\rank\alpha(x)-r
		.
	\]
	Now assume that also $(4)$ holds. Then, since $\frac{d\mu}{d\lambda}(x)=\infty$, by Proposition \ref{proposition 20}, (1), it follows that
	\[
		\lim_{z\downarrow x}\frac{\Im\tr M_w(z)}{\Im m(z)}=\frac{d\rho}{d\mu}(x)
		.
	\]
	Since the matrix $\lim_{z\downarrow x}\frac{\Im M_w(z)}{\Im m(z)}$ is positive semidefinite, we have 
	\begin{align*}
	\frac{d\rho}{d\mu}(x)=0\ \Leftrightarrow\ &\, 
	\tr\,\Big(\lim_{z\downarrow x}\frac{\Im M_w(z)}{\Im m(z)}\Big)=0\ \Leftrightarrow\ \lim_{z\downarrow x}\frac{\Im M_w(z)}{\Im m(z)}=0
	\\
	\Leftrightarrow\ &\, 
	\rank\,\Big(\lim_{z\downarrow x}\frac{\Im M_w(z)}{\Im m(z)}\Big)=0\ \Leftrightarrow\ \rank\alpha(x)=r.
	\end{align*}
	If $\rank\alpha(x)>r$, then $\frac{d\rho}{d\mu}(x)>0$ and it follows that the limit exists
	\begin{multline*}
		\lim_{z\downarrow x}\frac{\Im M_w(z)}{\Im\tr M_w(z)}
		=\lim_{z\downarrow x}\left(\frac{\Im M_w(z)}{\Im m(z)}\cdot\frac{\Im m(z)}{\Im\tr M_w(z)}\right)
		\\
		=
		\begin{pmatrix}
		\alpha_{11}(x)\!-\!\alpha_{11}(x)A_1
		\big(\alpha_{22}(x)\!+\!A_1^*\alpha_{11}(x)A_1\big)^{-1}A_1^*\alpha_{11}(x)
		& 0
		\\[2mm]
		0 & 0
		\end{pmatrix}
		\cdot\Big(\frac{d\rho}{d\mu}(x)\Big)^{-1}
	\end{multline*}
	and hence $\rank\,\lim_{z\downarrow x}\frac{\Im M_w(z)}{\Im\tr M_w(z)}=\alpha(x)-r$.
\end{proof}

\begin{remark}\label{rem r=n}
	Again consider the case $r=n$ separately. In this case $\alpha(x)=\alpha_{22}(x)$. If $\rank\alpha(x)=r=n$, then under conditions (1)--(3) 
	$$
	(A+BM_0(z))^{-1}=\frac1{m(z)}\left(\frac{A_{22}}{m(z)}+\frac{M_0(z)}{m(z)}\right)^{-1}\to0,\quad z\downarrow x,
	$$ 
	and owing to \eqref{ml2} $\frac{\Im M_w(z)}{\Im m(z)}\to0$ as $z\downarrow x$. If also (4) holds, then $\frac{d\rho}{d\mu}(x)=\lim_{z\downarrow x}\frac{\Im\tr M_w(z)}{\Im m(z)}=0$.
\end{remark}

\noindent
The proof of the many layers case of Theorem~\ref{theorem main}, {\rm(S2)--(S3)}, is now completed by observing that 
sufficiently many points $x\in\bb R$ satisfy conditions $(1)$--$(4)$ of Lemma~\ref{lemma 9}. 

\begin{lemma}\label{lemma 10}
	There exists a Borel set $W\subseteq\bb R$ with 
	\[
		\lambda(W)=\mu_s(\bb R\setminus W)=0
	\]
	such that for every $x\in W$
	\begin{enumerate}[$(1)$]
	\item $\frac{d\Omega_0}{d\mu}(x)$ exists,
	\item $\frac{d\mu}{d\lambda}(x)=\infty$, 
	\item $\lim\limits_{z\downarrow x}\frac{M_0(z)}{m(z)}$ exists and equals $\frac{d\Omega_0}{d\mu}(x)$,
	\item $\frac{d\rho}{d\mu}(x)$ exists and is finite, 
	\end{enumerate}
	and such that the functions $\frac{d\mu_l}{d\mu}$, $l=1,\ldots,n$, are measureable on $W$. 
\end{lemma}

\begin{proof}
	By Proposition \ref{prop DiB02} there exists a Borel set $X$ such that $\mu(\bb R\setminus X)=0$ and symmetric derivatives $\frac{d\mu_l}{d\mu}(x)$ exist for every $x\in X$, $l=1,\ldots,n$, and are measurable functions from $X$ to $[0,1]$ ($X$ is the intersection of such sets for each of the measures $\mu_l$). Then (1) holds for all $x\in X$.

	By Corollary \ref{prop Vallee-Poussin measures}, (1), items (1) and (4) hold $\mu$-a.e. and hence $\mu_s$-a.e. By
	Corollary \ref{prop Vallee-Poussin measures}, (4), (1), item (2) holds $\mu_s$-a.e. and on a set of zero Lebesgue
	measure. By Proposition \ref{proposition 20}, (3), item (3) holds $\mu_s$-a.e. Since the intersection of $\mu_s$-full sets is a $\mu_s$-full set, all of (1)--(4) hold $\mu_s$-a.e. and on a set of Lebesgue measure zero. There exists a Borel set $W$ contained in this intersection (including $X$) such that $\lambda(W)=0$ and $\mu_s(\bb R\setminus W)=0$. The functions $\frac{d\mu_l}{d\mu}$, $l=1,\ldots,n$, remain measurable on $W$.
\end{proof}

\begin{proof}[Proof of Theorem~\ref{theorem main}, {\rm(S2)--(S3)}]
	Consider the set $W$ from the above lemma. Define for $l=1,\ldots,n$ Borel sets
	\begin{equation}\label{equ 13}
	W^{(k)}:=\Big\{x\in W: \rank\alpha(x)=k\Big\},
	\end{equation}
	$$
	W^{(>k)}:=\bigcup_{l=k+1}^n W^{(l)},
	$$
	which are Borel owing to measurability of $\frac{d\mu_l}{d\mu}$ on $W$ for each $l$.
	
	Since $\{x\in\bb R:\exists\frac{d\Omega_0}{d\mu}(x),\rank\frac{d\Omega_0}{d\mu}(x)=k\}\setminus W^{(k)}$ are $\mu_s$-zero sets, we have $\mathds{1}_{\{N_0=k\}}\cdot\nu=\mathds{1}_{W^{(k)}}\cdot\nu$ for every measure $\nu\ll\mu_s$ (in particular, for $\rho_{s,ac}$).
	
	By Lemma~\ref{lemma 9} ( and Remark \ref{rem r=n}) we have $\frac{d\rho}{d\mu}(x)=0$ for all $x\in W^{(r)}$, hence by
	Corollary \ref{prop Vallee-Poussin sets}, (1), $\rho(W^{(r)})=0$. Therefore $\rho_{s,ac}(W^{(r)})=0$, which means that $\mathds{1}_{\{N_0=r\}}\cdot\rho_{s,ac}=0$,  and that is (S2).
	
	Further, by Lemma~\ref{lemma 9} we have $\frac{d\rho}{d\mu}(x)\in(0,\infty)$ for all $x\in W^{(>r)}$, and hence by
	Corollary \ref{prop Vallee-Poussin sets}, (4),
	$\mathds{1}_{W^{(>r)}}\cdot\mu\sim\mathds{1}_{W^{(>r)}}\cdot\rho$. Since $\lambda(W^{(>r)})=0$, it follows that
	$\mathds{1}_{W^{(>r)}}\cdot\mu_s=\mathds{1}_{W^{(>r)}}\cdot\mu\sim\mathds{1}_{W^{(>r)}}\cdot\rho=\mathds{1}_{W^{(>r)}}\cdot\rho_s$.
	As $\frac{d\rho}{d\mu}(x)\in[0,\infty)$ for all $x\in W$, taking into account (S1) which means that $\rho_{s,s}$ is the
	singular part in the Lebesgue decomposition of $\rho$ w.r.t. $\mu$, by Corollary \ref{prop Vallee-Poussin measures}, (4), we have $\rho_{s,s}(W)=0$. Then, in particular, $\mathds 1_{W^{(>r)}}\cdot\rho_{s,s}=0$. Therefore we have
	\begin{multline*}
	\mathds{1}_{\{N_0>r\}}\cdot\mu_s=\mathds{1}_{W^{(>r)}}\cdot\mu_s\sim\mathds{1}_{W^{(>r)}}\cdot\rho_s
	\\
	=\mathds{1}_{W^{(>r)}}\cdot\rho_{s,ac}+\mathds{1}_{W^{(>r)}}\cdot\rho_{s,s}=\mathds{1}_{W^{(>r)}}\cdot\rho_{s,ac}=\mathds{1}_{\{N_0>r\}}\cdot\rho_{s,ac},
	\end{multline*}
	which is the first part of (S3).
	
	For all $x\in W^{(>r)}$ (and therefore $\mathds{1}_{\{N_0>r\}}\cdot\rho_{s,ac}$-a.e.) by Lemma~\ref{lemma 9} it holds that 
	\[
	\rank\,\lim_{z\downarrow x}\frac{\Im M_w(z)}{\tr\Im M_w(z)}=\rank\alpha(x)-r=N_0(x)-r
	.
	\]
	By \eqref{equ 12} this means that
	\[
	N_{A,B}=N_0-r,\qquad \mathds{1}_{\{N_0>r\}}\cdot\rho_{s,ac}\text{-a.e.,}
	\]
	which completes the proof of {\rm (S3)}.
\end{proof}

\section{The few layers case}\label{section few layers}

In this section we prove assertions {\rm(S4)} and {\rm(S5)} from Theorem~\ref{theorem main}. 
This will be done by reducing to the already established many layers case {\rm(S2)--(S3)}, 
and referring  to the perturbation Lemma~\ref{lemma unknown}. 

The reduction is achieved by means of the following theorem, which is purely geometric. 
Recall that we denote for a pair $(A,B)$ of matrices which satisfies {\rm(D3)} the corresponding Lagrange plane 
as $\theta_{A,B}$, cf. \eqref{theta}.

\begin{theorem}\label{theorem 9}
	Let $A,B\in\bb C^{n\times n}$ such that $(A,B)$ satisfies {\rm(D3)} and {\rm(D4)}. 
	Then for each $k\in\{1,\ldots,\rank B\}$ there exist $A_k,B_k\in\bb C^{n\times n}$ such that 
	\begin{enumerate}[$(1)$]
	\item $(A_k,B_k)$ satisfies {\rm(D3)} and {\rm(D4)},
	\item $\rank B_k=k$ and 
		$\dim\left(\raisebox{2pt}{$\theta_{A,B}$}\big/\raisebox{-3pt}{$\theta_{A,B}\cap\theta_{A_k,B_k}$}\right)=\rank B-k$.
	\end{enumerate}
\end{theorem}

\noindent
The proof of this theorem relies on two lemmata. 

\begin{lemma}\label{lemma 7}
	Let $\theta$ and $\theta_0$ be Lagrange planes in $\bb C^n\times\bb C^n$, and let $\Pi$ be a linear subspace of $\bb C^n\times\bb C^n$ with 
	\[
		\theta\cap\theta_0\subseteq\Pi\subseteq\theta_0
		.
	\]
	Then there exists a Lagrange plane $\theta'$ in $\bb C^n\times\bb C^n$ such that 
	\[
		\theta'\cap\theta_0=\Pi,\quad 
		\dim\left(\raisebox{2pt}{$\theta$}\big/\raisebox{-3pt}{$\theta\cap\theta'$}\right)=\dim\Pi-\dim(\theta\cap\theta_0)
		.
	\]
\end{lemma}
\begin{proof}
	We have 
	\[
		\theta\cap\theta_0=\theta^{[\perp]}\cap\theta_0^{[\perp]}=(\theta+\theta_0)^{[\perp]}\subseteq
		(\theta+\Pi)^{[\perp]}
		.
	\]
	Choose a linear subspace $\Pi'$ with $(\theta+\Pi)^{[\perp]}=(\theta\cap\theta_0)\dot+\Pi'$. Observe that 
	\begin{equation}\label{equ 7}
		\Pi'\subseteq\theta^{[\perp]}=\theta,\quad \Pi'\cap\theta_0\subseteq\Pi'\cap\theta\cap\theta_0=\{0\}
		.
	\end{equation}
	In particular, we have $\Pi'\cap\Pi=\{0\}$. Now set $\theta':=\Pi'\dot+\Pi$. 
	Since $\Pi'$ and $\Pi$ are neutral subspaces with $\Pi'[\perp]\Pi$, the subspace $\theta'$ is neutral. 
	We have 
	\begin{equation}\label{equ 8}
		\theta\cap\Pi=\theta\cap\theta_0
		,
	\end{equation}
	and hence can compute 
	\begin{multline*}
		\dim(\theta\cap\theta_0)+\dim\Pi'=\dim(\theta+\Pi)^{[\perp]}=2n-\dim(\theta+\Pi)
		\\
		=2n-\big[\dim\theta+\dim\Pi-\dim(\theta\cap\Pi)\big]=n-\dim\Pi+\dim(\theta\cap\theta_0)
		.
	\end{multline*}
	Thus $\dim\theta'=\dim\Pi+\Pi'=n$, and we see that $\theta'$ is a Lagrange plane. 

	Using \eqref{equ 7} and\eqref{equ 8} we find 
	\begin{align*}
		& \theta'\cap\theta_0=(\Pi'\dot+\Pi)\cap\theta_0=(\Pi'\cap\theta_0)\dot+\Pi=\Pi
		,
		\\
		& \theta'\cap\theta=(\Pi'\dot+\Pi)\cap\theta=\Pi'\dot+(\theta_0\cap\theta)
		,
	\end{align*}
	and from the latter 
	\[
		\dim\left(\raisebox{2pt}{$\theta$}\big/\raisebox{-3pt}{$\theta\cap\theta'$}\right)=
		n-\big[\dim\Pi'+\dim(\theta_0\cap\theta)\big]=\dim\Pi-\dim(\theta_0\cap\theta)
		.
	\]
\end{proof}

\begin{lemma}\label{lemma 8}
	Let $\mc L$ be a linear subspace of $\bb C^n$ with $\dim\mc L<n$, and let $\mc C_1,\ldots,\mc C_N$ be a finite number of
	linear subspaces of $\bb C^n$ with 
	\[
		\forall j\in\{1,\ldots,N\}\colon \dim\mc C_j\leq n-\dim\mc L-1\text{ and }\mc L\cap\mc C_j=\{0\}
	\]
	Then there exists a linear subspace $\mc L'$ with 
	\[
		\mc L\subseteq\mc L',\quad \dim\mc L'=\dim\mc L+1,\quad \forall j\in\{1,\ldots,N\}\colon\mc L'\cap\mc C_j=\{0\}
	\]
\end{lemma}
\begin{proof}
	Choose a linear subspace $\mc K$ with $\bb C^n=\mc L\dot+\mc K$, and let $\pi:\bb C^n\to\mc K$ be the corresponding
	projection onto the second summand. 

	We first observe that for every subspace $\mc C$ with $\mc L\cap\mc C=\{0\}$ it holds that 
	\[
		\big\{z\in\mc K: (\mc L+\spn\{z\})\cap\mc C\neq\{0\}\big\}\subseteq\pi(\mc C)
		.
	\]
	Indeed, assume that $x\in\mc L,\lambda\in\bb C$, and $x+\lambda z\in\mc C\setminus\{0\}$. Since $\mc L\cap\mc C=\{0\}$, 
	we must have $\lambda\neq 0$, and hence $z=\frac 1\lambda\pi(x+\lambda z)\in\pi(\mc C)$. 

	Now let $\mc C_1,\ldots,\mc C_N$ be as in the statement of the lemma. Then 
	\[
		\Big\{z\in\mc K\ \big|\ 
		\exists j\in\{1,\ldots,N\}\colon(\mc L+\spn\{z\})\cap\mc C_j\neq\{0\}\Big\}
		\subseteq\bigcup_{j=1}^N\pi(\mc C_j)
		.
	\]
	Each $\pi(\mc C_j)$ is a linear subspace of $\mc K$ with 
	\[
		\dim\pi(\mc C_j)\leq\dim\mc C_j\leq n-\dim\mc L-1=\dim\mc K-1
		,
	\]
	and hence is a closed subset of $\mc K$ with empty interior. Thus also the finite union $\bigcup_{j=1}^N\pi(\mc C_j)$ 
	is a subset of $\mc K$ with empty interior, in particular, 
	\[
		\bigcup_{j=1}^N\pi(\mc C_j)\neq\mc K
		.
	\]
	For each element $z\in\mc K\setminus\bigcup_{j=1}^N\pi(\mc C_j)$, the space $\mc L':=\mc L+\spn\{z\}$ has the required
	properties. 
\end{proof}

\begin{proof}[Proof of Theorem~\ref{theorem 9}]
	We use induction on $k$. For $k=\rank B$ there is nothing to prove: just set $A_0:=A$ and $B_0:=B$.
	Assume we have $k\in\{1,\ldots,\rank B\}$ and $(A_k,B_k)$ with $(1)$ and $(2)$. The aim is to construct
	$(A_{k-1},B_{k-1})$. 

	We work with Lagrange planes rather than matrices: denote 
	\[
		\theta:=\theta_{A,B},\quad \theta_k:=\theta_{A_k,B_k},\quad \theta_0=\{0\}\times\bb C^n
		.
	\]
	Further, set $\mc L:=\ker B_k$, so that $\theta_k\cap\theta_0=\{0\}\times\mc L$. Then $\dim\mc L=n-k<n$, and by {\rm(D4)} 
	we have $\mc L\cap\mc C=\{0\}$ for every coordinate plane $\mc C$ with dimension $k-1$. According to 
	Lemma~\ref{lemma 8} we find $\mc L'\supseteq\mc L$ with $\dim\mc L'=n-k+1$ such that still $\mc L'\cap\mc C=\{0\}$ 
	for all coordinate planes $\mc C$ with dimension $k-1$. 

	Now set $\Pi:=\{0\}\times\mc L'$. According to Lemma~\ref{lemma 7} we find a Lagrange plane $\theta'$ with 
	\begin{equation}\label{equ 9}
		\theta'\cap\theta_0=\{0\}\times\mc L',\quad
		\dim\left(\raisebox{2pt}{$\theta_k$}\big/\raisebox{-3pt}{$\theta_k\cap\theta'$}\right)=
		\dim\Pi-\dim(\theta_k\cap\theta_0)=1
		.
	\end{equation}
	We need to compute the dimension of the factor space $\theta/(\theta\cap\theta')$. The canonical map 
	\[
		\raisebox{2pt}{$\theta'\cap\theta_0$}\big/\raisebox{-3pt}{$(\theta'\cap\theta)\cap\theta_0$}
		\ \longrightarrow\ 
		\raisebox{2pt}{$\theta'$}\big/\raisebox{-3pt}{$\theta'\cap\theta$}
	\]
	is injective, and $\dim\theta=\dim\theta'$. Hence, 
	\begin{align*}
		\dim\left(\raisebox{2pt}{$\theta$}\big/\raisebox{-3pt}{$\theta\cap\theta'$}\right)
		=&\, \dim\left(\raisebox{2pt}{$\theta'$}\big/\raisebox{-3pt}{$\theta\cap\theta'$}\right)
		\geq\dim\left(\raisebox{2pt}{$\theta'\cap\theta_0$}\big/\raisebox{-3pt}{$\theta\cap\theta'\cap\theta_0$}\right)
		\\
		= &\, \dim(\theta'\cap\theta_0)-\dim(\theta\cap\theta'\cap\theta_0)
		\geq\dim(\theta'\cap\theta_0)-\dim(\theta\cap\theta_0)
		\\
		=&\, \dim\mc L'-(n-\rank B)=\rank B-k+1
		.
	\end{align*}
	On the other hand, using \eqref{equ 9}, $\dim\theta_k=\dim\theta'$, and the respective canonical injections, yields
	\[
		\dim\left(\raisebox{2pt}{$\theta_k\cap\theta$}\big/\raisebox{-3pt}{$\theta'\cap\theta_k\cap\theta$}\right)\leq 1
		,\quad 
		\dim\left(\raisebox{2pt}{$\theta'\cap\theta$}\big/\raisebox{-3pt}{$\theta'\cap\theta_k\cap\theta$}\right)\leq 1
		.
	\]
	From 
	\begin{multline*}
		\left(\raisebox{2pt}{$\theta$}\big/\raisebox{-3pt}{$\theta\cap\theta_k$}\right)\times
		\left(\raisebox{2pt}{$\theta\cap\theta_k$}\big/\raisebox{-3pt}{$\theta\cap\theta_k\cap\theta'$}\right)
		\\
		\cong\,
		\left(\raisebox{2pt}{$\theta$}\big/\raisebox{-3pt}{$\theta\cap\theta_k\cap\theta'$}\right)
		\\
		\cong\,
		\left(\raisebox{2pt}{$\theta$}\big/\raisebox{-3pt}{$\theta\cap\theta'$}\right)\times
		\left(\raisebox{2pt}{$\theta\cap\theta'$}\big/\raisebox{-3pt}{$\theta\cap\theta_k\cap\theta'$}\right)
	\end{multline*}
	we now obtain 
	\[
		\Big|\dim\left(\raisebox{2pt}{$\theta$}\big/\raisebox{-3pt}{$\theta\cap\theta_k$}\right)-
		\dim\left(\raisebox{2pt}{$\theta$}\big/\raisebox{-3pt}{$\theta\cap\theta'$}\right)\Big|\leq 1
		,
	\]
	and hence 
	\[
		\dim\left(\raisebox{2pt}{$\theta$}\big/\raisebox{-3pt}{$\theta\cap\theta'$}\right)
		\leq
		\dim\left(\raisebox{2pt}{$\theta$}\big/\raisebox{-3pt}{$\theta\cap\theta_k$}\right)+1\leq \rank B-k+1
		.
	\]
	It remains to choose $(A_{k-1},B_{k-1})$ with $\theta'=\theta_{A_{k-1},B_{k-1}}$. Then $(A_{k-1},B_{k-1})$ satisfies 
	$(1)$ and $(2)$. 
\end{proof}

\noindent
It is now easy to deduce the few layers case in Theorem~\ref{theorem main}.

\begin{proof}[Proof of Theorem~\ref{theorem main}, {\rm(S1), (S4)--(S5)}]
	First note that (S1) follows from \cite{Gesztesy-Tsekanovskii-2000}, as we mentioned above, and (S5) is obtained by application of Lemma \ref{lemma unknown}.
	
	 Let us prove {\rm(S1)}. Consider $k\in\{1,\ldots,r-1\}$. 
	Choose $(A_k,B_k)$ according to Theorem~\ref{theorem 9}. 
	Denote by $\rho_k$ the scalar spectral measure of $L_{A_k,B_k}$ (i.e., $\rho_k=\tr\Omega_{A_k,B_k}$), let $W_k\subseteq\bb R$ be a set as constructed in
	Lemma~\ref{lemma 10} for $(A_k,B_k)$, and let $W^{(k)}_k$ be the corresponding set from \eqref{equ 13}.
	Then we know, from the already proven part {\rm(S2)} that $\rho_k(W^{(k)}_k)=0$. 

	Let $\rho=\rho_{ac,k}+\rho_{s,k}$ be the Lebesgue decomposition of $\rho$ w.r.t.\ $\rho_k$. Since $\rho_k(W^{(k)}_k)=0$,
	it follows that $\rho_{ac,k}(W^{(k)}_k)=0$ or, equivalently, $\mathds 1_{W^{(k)}_k}\cdot\rho_{ac,k}=0$, hence
	\begin{equation}\label{200}
		\mathds{1}_{W^{(k)}_k}\cdot\rho_{s,ac}\ll\mathds{1}_{W^{(k)}_k}\cdot\rho=\mathds 1_{W^{(k)}_k}\cdot\rho_{ac,k}+\mathds 1_{W^{(k)}_k}\cdot\rho_{s,k}=\mathds 1_{W^{(k)}_k}\cdot\rho_{s,k}\ll\rho_{s,k}
		.
	\end{equation}
	Lemma \ref{lemma ranks} applied to $A,B$ and $A_k,B_k$ gives 
	$$
	\rank[(L_{A,B}-\lambda I)^{-1}-(L_{A_k,B_k}-\lambda I)^{-1}]\leq\dim\left(\raisebox{2pt}{$\theta_{A,B}$}\big/\raisebox{-3pt}{$\theta_{A,B}\cap\theta_{A_k,B_k}$}\right)=r-k
	$$ 
	and Lemma~\ref{lemma unknown} applied to $L_{A,B}$ and $L_{A_k,B_k}$ gives
	\[
		N_{A,B}\leq r-k,\qquad\text{$\rho_{s,k}$-a.e.}
	\]
	and, owing to \eqref{200}, $\mathds{1}_{W^{(k)}_k}\cdot\rho_{s,ac}$-a.e. On the other hand, from the definitions of sets $W^{(k)}_k$ we have 
	\[
		N_0=k,\qquad\text{$\mathds{1}_{W^{(k)}_k}\cdot\mu$-a.e.\ and hence $\mathds{1}_{W^{(k)}_k}\cdot\rho_{s,ac}$-a.e.,}
	\]
	therefore
	\begin{equation}\label{100}
		N_{A,B}\leq r-N_0,\qquad\text{$\mathds{1}_{(\cup_{k=1}^{r-1}W^{(k)}_k)}\cdot\rho_{s,ac}$-a.e.}
	\end{equation}
	Consider the intersection 
	$$
	\widetilde W:=\bigcap_{k=1}^{r-1}W_k
	$$ 
	which is also a $\lambda$-zero and $\mu_s$-full set, and analogously to \eqref{equ 13} define $\widetilde W^{(k)}:=\{x\in\widetilde W:\rank\alpha(x)=k\}$ for $k=1,\ldots,n$. Then we have $\mathds{1}_{(\cup_{k=1}^{r-1}\widetilde W^{(k)})}\cdot\mu_s=\mathds{1}_{\{0<N_0<r\}}\cdot\mu_s$. Since for every $k=1,\ldots,n$ one has $\mu_s(W_k^{(k)}\setminus\widetilde W^{(k)})=0$, it is true that $\mathds 1_{(\cup_{k=1}^{r-1}W_k^{(k)})}\cdot\mu_s=\mathds 1_{\{0<N_0<r\}}\cdot\mu_s$ as well, and therefore also $\mathds 1_{(\cup_{k=1}^{r-1}W_k^{(k)})}\cdot\rho_{s,ac}=\mathds 1_{\{0<N_0<r\}}\cdot\rho_{s,ac}$. Thus \eqref{100} in fact coincides with the assertion of (S5), which completes the proof.
\end{proof}

\section*{Declarations}
There are no conflicts of interest.
Data sharing not applicable to this article as no datasets were generated or analysed during the current study.

{\footnotesize
	\begin{flushleft}
		S.\,Simonov\\
		St. Petersburg Department of V. A. Steklov Institute of Mathematics of the Russian Academy of Sciences, 
		Fontanka 27, St. Petersburg 191023\\
		and St. Petersburg State University, 
		Universitetskaya nab. 7--9, St. Petersburg 199034\\
		RUSSIA\\
		email: sergey.a.simonov@gmail.com
	\end{flushleft}
	\begin{flushleft}
		H.\,Woracek\\
		Institut f\"ur Analysis und Scientific Computing\\
		Technische Universit\"at Wien\\
		Wiedner Hauptstra{\ss}e.\ 8--10/101\\
		1040 Wien\\
		AUSTRIA\\
		email: harald.woracek@tuwien.ac.at
	\end{flushleft}

}

\end{document}